\documentclass[11pt,letterpaper]{article}
\usepackage[margin=1in]{geometry}

\usepackage[T1]{fontenc}
\usepackage{lmodern}
\usepackage{microtype}

\usepackage{amsmath,amssymb,amsthm,amsfonts,mathtools}
\usepackage{float}
\usepackage{placeins}

\usepackage{graphicx}
\graphicspath{{figures/}}

\usepackage{subcaption}
\usepackage{algorithm}
\usepackage{algpseudocode}

\usepackage{booktabs,longtable}
\usepackage[numbers]{natbib}
\usepackage{xurl}
\usepackage[hidelinks]{hyperref}
\usepackage{doi}
\newtheorem{theorem}{Theorem}[section]

\newtheorem{lemma}[theorem]{Lemma}

\newtheorem{proposition}[theorem]{Proposition}

\newtheorem{corollary}[theorem]{Corollary}

\newtheorem{assumption}[theorem]{Assumption}

\newtheorem{definition}[theorem]{Definition}

\newcommand{\R}{\mathbb{R}}

\newcommand{\E}{\mathbb{E}}
\newcommand{\Prob}{\mathbb{P}}

\newcommand{\1}{\mathbf{1}}

\DeclareMathOperator{\Tr}{Tr}

\DeclareMathOperator{\diag}{diag}

\DeclareMathOperator*{\argmin}{arg\,min}
\DeclareMathOperator*{\argmax}{arg\,max}

\newcommand{\eps}{\varepsilon}

\title{Certified Approximation of Probabilistic Reachable Sets via Parabolic Hamilton--Jacobi--Isaacs Equations}
\author{%
  Minjung Gim\textsuperscript{1}\thanks{Minjung Gim and Hee Jun Yang contributed equally to this work.}%
  \qquad Hee Jun Yang\textsuperscript{1}\footnotemark[1]%
  \qquad Yeoneung Kim\textsuperscript{2}\thanks{Corresponding author.
    \href{mailto:yeoneung@yonsei.ac.kr}{\nolinkurl{yeoneung@yonsei.ac.kr}}}%
  \\[0.6em]
  \small\textsuperscript{1}National Institute for Mathematical Sciences (NIMS)\\
  \small 70, Yuseong-daero 1689 beon-gil, Yuseong-gu, Daejeon 34047, Republic of Korea\\[0.25em]
  \small\textsuperscript{2}Department of Industrial Engineering, Yonsei University\\
  \small 50 Yonsei-ro, Seodaemun-gu, Seoul 03722, Republic of Korea%
}
\date{}

\hypersetup{
  pdftitle={Certified Approximation of Probabilistic Reachable Sets via Parabolic Hamilton--Jacobi--Isaacs Equations},
  pdfauthor={Minjung Gim, Hee Jun Yang, Yeoneung Kim},
  pdfkeywords={Probabilistic reachability, Hamilton--Jacobi--Isaacs equations, Stochastic differential games, A posteriori certification, Physics-informed neural networks, Policy iteration}
}

\begin{document}
\maketitle

\begin{abstract}
We study terminal-time probabilistic reachability for controlled diffusions with adversarial disturbances. Our main result converts verified one-sided residual bounds for the full parabolic Hamilton--Jacobi--Isaacs operator into certified inner and outer reachable-set approximations. One-sided terminal envelopes handle the discontinuous target payoff, while explicit lateral inequalities make the certificate computable on a bounded domain. Neural policy iteration generates candidates; certification is independent of its empirical training loss. The experiments serve two complementary objectives. A scalar game demonstrates complete outward-rounded interval certification and a nonempty certified inner set. Air3D and quadrotor benchmarks demonstrate neural candidate generation in three, four, and ten state dimensions through reference comparisons, feedback rollouts, and adaptive sampling. On the evaluated 10D slices, adaptive sampling achieves approximately one-percent relative error against the restricted-product reference. Exact one-sided Howard convergence under explicit regularity and compactness assumptions supports candidate generation.
\end{abstract}


\section{Introduction}
Reachability analysis for systems subject to random disturbances asks how
likely the state is to enter a target. When an adversary also influences
the dynamics, this probability depends on competing control decisions.
We consider a terminal-time problem: a controller seeks to avoid a target
at a prescribed time, while a disturbance seeks to force the terminal
state into it. The aim is to identify initial states from which the
disturbance can guarantee a specified probability of reaching the target
against every admissible control.

Stochastic reachability and reach-avoid analysis have established
probability-based verification methods for controlled and hybrid systems
\cite{abate2008probabilistic,SummersLygeros2010,EsfahaniChatterjeeLygeros2016}.
Related approaches combine deterministic reachability with stochastic
deviation bounds \cite{Jafarpour2025PRS}, or use stochastic barrier functions
\cite{SantoyoDutreixCoogan2021}. In a stochastic differential game, the
terminal-event probability is represented by a value function governed by
a parabolic Hamilton--Jacobi--Isaacs (HJI) equation
\cite{FlemingSouganidis1989,FlemingSoner2006,BuckdahnLi2008}.
The target indicator gives discontinuous terminal data, for which
viscosity theory provides relaxed terminal traces and comparison tools
\cite{CrandallIshiiLions1992,BarlesPerthame1987,Barles1994}.

Neural PDE solvers, including physics-informed neural networks (PINNs)
\cite{RaissiPerdikarisKarniadakis2019}, deep Galerkin methods
\cite{SirignanoSpiliopoulos2018}, and DeepReach \cite{bansal2021deepreach},
provide candidate value functions beyond the dimensions accessible to dense
grids. For nonconvex stochastic HJI equations, Yang et al.\
\cite{YangGimKim2027NonconvexHJI} combine neural evaluation of fixed feedback
pairs with pointwise minimax policy updates. Yang and Kim
\cite{YangKim2026FiniteDifferencePI} use monotone finite differences for
neural policy iteration on degenerate Hamilton--Jacobi--Bellman (HJB) and HJI
equations.

Turning the level sets of neural candidates into probability guarantees
requires verified error bounds. Residual verification has been developed for deterministic
reachability and optimal control: Solanki et al.\ \cite{Solanki2025CARe}
use SMT-based residual bounds and counterexamples, and Liu
\cite{Liu2026Verifiable} derives verifiable residual error bounds for
Lyapunov and HJB equations. Gomes et al.\
\cite{GomesGuoKimTranZhang2026NewtonPI} use residual estimates to control the
error of PINN approximations to individual Newton corrections for viscous
Hamilton--Jacobi equations. For the stochastic game considered here,
verification must also account for diffusion, the discontinuous target
payoff, and a reachable set defined through nonanticipative disturbance
strategies.

We establish an a posteriori certificate that converts verified one-sided
residual bounds for the full parabolic HJI operator into inner and outer
reachable-set approximations. One-sided terminal envelopes bracket the
target indicator, and explicit lateral inequalities permit verification
on a bounded domain. The certificate applies to any sufficiently regular
candidates satisfying these conditions, independently of their training
procedure or empirical loss.

We generate candidates by physics-informed neural network policy iteration
(PINN--PI) and prove convergence
of its exact one-sided Howard counterpart under explicit regularity and
compactness assumptions. The experiments address certification and
candidate generation separately.
A scalar game demonstrates the complete outward-rounded interval
certificate with a nonempty inner set. Air3D and quadrotor experiments
evaluate candidate generation in three, four, and ten state dimensions
through reference comparisons, feedback rollouts, and adaptive sampling.

\section{Reachability game and terminal envelopes}
\label{sec:formulation}
Continuous terminal envelopes connect the target event to value functions
that can be compared with smooth candidates.

\subsection{Controlled diffusion and reachability}
Let \(W\) be an \(m\)-dimensional Brownian motion on canonical Wiener
space. For a game starting at \(t\), use the usual augmentation of the
filtration generated by \((W_s-W_t)_{s\ge t}\), and consider
\begin{equation}\label{eq:SDE}
 dX_s=f(s,X_s,a_s,b_s)\,ds+\sigma(s,X_s)\,dW_s,
 \qquad X_t=x,\quad t\le s\le T.
\end{equation}
\begin{assumption}\label{ass:dynamics}
The action sets \(A\) and \(B\) are nonempty compact metric spaces.
The drift \(f:[0,T]\times\R^d\times A\times B\to\R^d\) and diffusion
\(\sigma:[0,T]\times\R^d\to\R^{d\times m}\) are continuous, globally
Lipschitz in \(x\) uniformly in the other variables, and of linear growth.
The covariance \(\Sigma=\sigma\sigma^\top\) is uniformly elliptic:
\(z^\top\Sigma(t,x)z\ge\lambda|z|^2\) for some \(\lambda>0\) and all
\((t,x,z)\).
\end{assumption}
These assumptions give a unique strong solution \(X^{t,x;a,b}\) for each
admissible pair. Let \(\mathcal A_{t,s}\) and \(\mathcal B_{t,s}\) be the
progressively measurable \(A\)- and \(B\)-valued processes on \([t,s]\).
A disturbance strategy \(\gamma:\mathcal A_{t,s}\to\mathcal B_{t,s}\)
is nonanticipative if agreement of two controls on \([t,\tau]\) implies
agreement of their images on \([t,\tau]\), for every stopping time
\(\tau\) taking values in \([t,s]\).
Process equalities are understood \(dt\otimes d\mathbb P\)-a.e.
Write \(\Gamma_{t,s}\) for these strategies, and use the subscript \(t\)
when \(s=T\). The admissible classes are stable under restriction and
concatenation at deterministic times.

For a bounded Borel payoff \(\psi\), define
\begin{equation}\label{eq:general_game_value}
 V^\psi(t,x)=\sup_{\gamma\in\Gamma_t}\inf_{a\in\mathcal A_t}
 \E[\psi(X_T^{t,x;a,\gamma[a]})].
\end{equation}
Let \(\mathcal G\subsetneq\R^d\) be nonempty and closed, set
\(g=\1_{\mathcal G}\), and write
\begin{equation}\label{eq:indicator_game_value}
 v=V^g,\qquad 0\le v\le1.
\end{equation}
The controller minimizes the terminal hitting probability and the
disturbance maximizes it. For \(\delta\in(0,1]\), the relevant set is
\begin{definition}\label{def:deltaBRS}
\[
 R_\delta(t)=\{x:\exists\gamma\in\Gamma_t\ \forall a\in\mathcal A_t,
 \ \Prob(X_T^{t,x;a,\gamma[a]}\in\mathcal G)\ge1-\delta\}.
\]
\end{definition}
This is a terminal event; paths need not avoid or enter the target at
intermediate times.

\begin{proposition}\label{prop:value_level_reachability}
For every \(\xi>0\),
\[
 \{x:v(t,x)\ge1-\delta+\xi\}\subset R_\delta(t)
 \subset\{x:v(t,x)\ge1-\delta\}.
\]
If the supremum defining \(v\) is attained, the right-hand inclusion is
an equality and the inner margin may be omitted.
\end{proposition}
\begin{proof}
A strategy guaranteeing probability at least \(1-\delta\) implies the
same lower bound on its supremum. Conversely, if
\(v(t,x)\ge1-\delta+\xi\), choose a strategy whose infimum payoff is
greater than \(v(t,x)-\xi\). It guarantees the required probability.
Attainment permits \(\xi=0\).
\end{proof}

\subsection{HJI equation and terminal envelopes}
The local minimax drift is encoded by the Hamiltonian
\begin{equation}\label{eq:Isaacs_Hamiltonian}
 H(t,x,p)=\sup_{b\in B}\inf_{a\in A} f(t,x,a,b)\cdot p.
\end{equation}
\begin{assumption}\label{ass:isaacs}
For every \((t,x,p)\),
\(\sup_b\inf_a f(t,x,a,b)\cdot p=\inf_a\sup_b f(t,x,a,b)\cdot p\).
\end{assumption}
We use the backward operator
\[
 \mathcal H[u]=\mathcal L u
 :=\partial_tu+H(t,x,\nabla_xu)+\tfrac12\Tr(\Sigma D^2_{xx}u).
\]
Assumptions~\ref{ass:dynamics} and~\ref{ass:isaacs} remain in force
throughout the theoretical results.
Thus a classical subsolution satisfies \(\mathcal H[u]\ge0\), and a
supersolution satisfies \(\mathcal H[u]\le0\). For bounded uniformly
continuous \(\psi\), the standard game representation gives
the viscosity solution \(U=V^\psi\) of
\begin{equation}\label{eq:HJI_general_payoff}
 \mathcal H[U]=0\quad\text{on }[0,T)\times\R^d,
 \qquad U(T,\cdot)=\psi.
\end{equation}
The supporting comparison and terminal-trace results are recorded in
Appendix~\ref{sec:HJI}.

The indicator payoff is handled by continuous envelopes. Let
\(\rho(x)=\operatorname{dist}(x,\mathcal G^c)
-\operatorname{dist}(x,\mathcal G)\), and let
\(\chi\in C^\infty(\R;[0,1])\) be nondecreasing, with \(\chi=0\) on
\(( -\infty,0]\) and \(\chi=1\) on \([1,\infty)\). Set
\begin{equation}\label{eq:one_sided_terminal_envelopes}
 g^-_\eps(x)=\chi(\rho(x)/\eps),\qquad
 g^+_\eps(x)=\chi((\rho(x)+\eps)/\eps).
\end{equation}
These Lipschitz functions satisfy
\begin{equation}\label{eq:terminal_ordering}
 g^-_\eps\le\1_{\operatorname{int}\mathcal G}
 \le\1_{\mathcal G}\le g^+_\eps,
\end{equation}
and differ only where \(|\rho|<\eps\). Smooth replacements may be used
provided the same one-sided ordering is preserved. Define
\(v^\pm_\eps=V^{g^\pm_\eps}\).
Expectation, infimum, and supremum preserve the payoff ordering
\eqref{eq:terminal_ordering}, so
\begin{equation}\label{eq:value_bracketing}
 v^-_\eps(t,x)\le v(t,x)\le v^+_\eps(t,x).
\end{equation}
Under the boundary-tube condition below,
Corollary~\ref{cor:v_relaxed_solution} establishes convergence of these
envelopes and uniqueness for the indicator-terminal HJI problem.

\subsection{Envelope gap}

To quantify the envelope gap, define
\[
    \mathcal N_\eps(\partial\mathcal G)
    :=
    \{x\in\R^d:|\rho(x)|<\eps\}.
\]

\begin{assumption}
\label{ass:tube_estimate}
For every compact set \(K\subset\R^d\) and every \(\tau\in(0,T]\), there
exists a modulus \(\omega_{K,\tau}(\eps)\to0\) as \(\eps\downarrow0\)
such that
\[
    \sup_{\substack{(t,x)\in[0,T-\tau]\times K\\
                    a\in\mathcal A_t,\ \gamma\in\Gamma_t}}
    \Prob\!\left(
        X_T^{t,x;a,\gamma[a]}
        \in
        \mathcal N_\eps(\partial\mathcal G)
    \right)
    \le
    \omega_{K,\tau}(\eps).
\]
\end{assumption}

\begin{proposition}
\label{prop:envelope_error}
Under Assumption~\ref{ass:tube_estimate}, for every compact
\(K\subset\R^d\), every \(\tau\in(0,T]\), and every
\((t,x)\in[0,T-\tau]\times K\),
\begin{equation}
\label{eq:envelope_error}
    0
    \le
    v^+_\eps(t,x)-v^-_\eps(t,x)
    \le
    \omega_{K,\tau}(\eps).
\end{equation}
In particular, \(v^-_\eps\) and \(v^+_\eps\) converge to the same limit
locally uniformly away from the terminal time whenever
\(\omega_{K,\tau}(\eps)\to0\).
\end{proposition}

The proof is given in Appendix~\ref{sec:HJI}.

If the terminal laws admit density upper bounds uniform over the admissible
controls and strategies for \((t,x)\in[0,T-\tau]\times K\), then a tube
volume tending to zero implies Assumption~\ref{ass:tube_estimate}.
For bounded \(C^2\) targets, such a uniform density bound gives
\(\omega_{K,\tau}(\eps)\le C_{K,\tau}\eps\). The required uniformity
over the admissible class is part of this sufficient condition.

The vanishing gap also identifies the indicator value through viscosity
stability and comparison.

\begin{corollary}
\label{cor:v_relaxed_solution}
Assume Assumptions~\ref{ass:dynamics}, \ref{ass:isaacs},
and~\ref{ass:tube_estimate}.  Then \(v^\pm_\eps\to v\) locally uniformly
on \([0,T)\times\R^d\) as \(\eps\downarrow0\); in particular, \(v\) is
continuous on \([0,T)\times\R^d\).  Moreover, \(v\) is the unique bounded relaxed
viscosity solution of \eqref{eq:HJI} in the sense of
Definition~\ref{def:relaxed_viscosity_solution}.
\end{corollary}

\begin{proof}
The ordering \eqref{eq:value_bracketing} and
Proposition~\ref{prop:envelope_error} give \(\|v^\pm_\eps-v\|_{L^\infty([0,T-\tau]\times K)}
\le\omega_{K,\tau}(\eps)\to0\) for every compact \(K\) and
\(\tau\in(0,T]\). The continuous envelope solutions therefore converge
locally uniformly to a continuous \(v\). Viscosity stability
\cite{CrandallIshiiLions1992} gives \(\mathcal Lv=0\) in the interior;
Theorem~\ref{thm:relaxed_wellposedness}(i) supplies the relaxed terminal
inequalities. Thus \(v\) is a bounded relaxed solution of \eqref{eq:HJI}.

For uniqueness, let \(w\) be any bounded relaxed viscosity solution.
The terminal ordering \eqref{eq:terminal_ordering} and
Definition~\ref{def:relaxed_viscosity_solution} give
\[
 g^-_\eps\le g_*\le w_*(T,\cdot),
 \qquad w^*(T,\cdot)\le g^*\le g^+_\eps.
\]
Since \(v^\pm_\eps\) attain \(g^\pm_\eps\) continuously at \(T\),
Theorem~\ref{thm:relaxed_wellposedness}(ii) applied to
\((v^-_\eps,w)\) and \((w,v^+_\eps)\) yields
\[
 v^-_\eps\le w\le v^+_\eps
 \qquad\text{on }[0,T)\times\R^d.
\]
Both \(w\) and \(v\) lie in this bracket, so
\(\sup_{[0,T-\tau]\times K}|w-v|
\le\omega_{K,\tau}(\eps)\).
Letting \(\eps\downarrow0\) gives \(w=v\) throughout
\([0,T)\times\R^d\).
\end{proof}

\section{A posteriori reachable-set certificates}\label{sec:certificate}
\label{subsec:one_sided_aposteriori}

The certificate converts residual and boundary inequalities into bounds
on the reachability value, independently of how candidates are generated.

\subsection{Whole-space certificate}
One-sided residual errors are absorbed by corrections proportional to
the remaining time.

\begin{theorem}
\label{thm:one_sided_certificate}
Assume Assumptions~\ref{ass:dynamics} and~\ref{ass:isaacs}, so that the
comparison principle of Theorem~\ref{thm:relaxed_wellposedness}(ii)
applies.  Let
\[
    \underline u,\overline u
    \in
    C([0,T]\times\R^d)
    \cap
    C^{1,2}([0,T)\times\R^d)
\]
be bounded functions satisfying
\[
    \underline u(T,x)\le g^-_\eps(x),
    \qquad
    \overline u(T,x)\ge g^+_\eps(x),
    \qquad x\in\R^d.
\]
Suppose that, for some constants \(\eta_-,\eta_+\ge0\),
\[
    \mathcal H[\underline u](t,x)\ge -\eta_-,
    \qquad
    \mathcal H[\overline u](t,x)\le \eta_+
\]
for all \((t,x)\in[0,T)\times\R^d\).  Then
\begin{equation}
\label{eq:certified_value_bounds}
    \underline u(t,x)-\eta_-(T-t)
    \le
    v^-_\eps(t,x)
    \le
    v(t,x)
    \le
    v^+_\eps(t,x)
    \le
    \overline u(t,x)+\eta_+(T-t)
\end{equation}
for all \((t,x)\in[0,T]\times\R^d\).

For every \(t\in[0,T]\), \(\delta\in(0,1]\), and \(\xi>0\),
\[\begin{gathered}
    \left\{
        x:
        \underline u(t,x)-\eta_-(T-t)
        \ge
        1-\delta+\xi
    \right\}
    \subset
    R_\delta(t)
    \\
    \subset
    \left\{
        x:
        \overline u(t,x)+\eta_+(T-t)
        \ge
        1-\delta
    \right\}.
\end{gathered}\]
If the supremum in the definition of \(v(t,x)\) is attained, the margin
\(\xi>0\) in the inner approximation may be removed.
\end{theorem}

\begin{proof}
Set \(\underline U=\underline u-\eta_-(T-t)\) and
\(\overline U=\overline u+\eta_+(T-t)\). Since the corrections are spatially
constant,
\[
 \mathcal H[\underline U]=\mathcal H[\underline u]+\eta_-\ge0,
 \qquad
 \mathcal H[\overline U]=\mathcal H[\overline u]-\eta_+\le0.
\]
Their terminal values are ordered below \(g^-_\eps\) and above
\(g^+_\eps\), respectively. Comparison with the continuous envelope
solutions gives \(\underline U\le v^-_\eps\) and
\(v^+_\eps\le\overline U\). The ordering \eqref{eq:value_bracketing}
supplies the middle inequalities. Applying
Proposition~\ref{prop:value_level_reachability} to these value bounds
gives the set inclusions and the attainment qualification.
\end{proof}

This certificate and its bounded-domain version in
Corollary~\ref{cor:localized_certificate} use continuous-envelope comparison
and payoff ordering, without Assumption~\ref{ass:tube_estimate}.
That assumption is used to make the envelope gap vanish and establish
uniqueness for the indicator-terminal problem in
Corollary~\ref{cor:v_relaxed_solution}.

\subsection{Bounded-domain certificate}
Verification is performed on a finite cylinder. Since the envelope values
lie in \([0,1]\), explicit lateral inequalities provide the ordering needed
to compare with their restrictions to that cylinder.
\begin{corollary}
\label{cor:localized_certificate}
Let \(\Omega\subset\R^d\) be a bounded domain and
\(Q:=[0,T)\times\Omega\).  Let
\(\underline u,\overline u\in C(\bar Q)\cap C^{1,2}(Q)\) satisfy, for
some constants \(\eta_-,\eta_+\ge0\), the interior residual
inequalities
\[
    \mathcal H[\underline u]\ge-\eta_-,
    \qquad
    \mathcal H[\overline u]\le\eta_+
    \qquad\text{on }Q,
\]
the terminal inequalities
\[
    \underline u(T,x)\le g^-_\eps(x),
    \qquad
    \overline u(T,x)\ge g^+_\eps(x),
    \qquad x\in\Omega,
\]
and the lateral boundary inequalities
\[
    \underline u(t,x)\le 0,
    \qquad
    \overline u(t,x)\ge 1,
    \qquad
    (t,x)\in[0,T)\times\partial\Omega .
\]
Then
\[\begin{gathered}
    \underline u(t,x)-\eta_-(T-t)
    \le
    v^-_\eps(t,x)
    \le
    v(t,x)
    \\
    \le
    v^+_\eps(t,x)
    \le
    \overline u(t,x)+\eta_+(T-t)
    \qquad\text{on }Q.
\end{gathered}\]
Consequently, for every \(t\in[0,T)\), \(\delta\in(0,1]\), and
\(\xi>0\),
\[
    \bigl\{x\in\Omega:
    \underline u(t,x)-\eta_-(T-t)\ge1-\delta+\xi\bigr\}
    \subset
    R_\delta(t),
\]
\[
    R_\delta(t)\cap\Omega
    \subset
    \bigl\{x\in\Omega:
    \overline u(t,x)+\eta_+(T-t)\ge1-\delta\bigr\}.
\]
\end{corollary}

\begin{proof}
Use the corrected candidates \(\underline U=\underline u-\eta_-(T-t)\)
and \(\overline U=\overline u+\eta_+(T-t)\) from
Theorem~\ref{thm:one_sided_certificate}. They are a subsolution and a
supersolution on \(Q\), with the required terminal ordering. Since
\(0\le v^\pm_\eps\le1\), the lateral conditions give
\[
 \underline U\le0\le v^-_\eps,
 \qquad v^+_\eps\le1\le\overline U
 \quad\text{on }[0,T)\times\partial\Omega.
\]
The envelope solutions are continuous by
Theorem~\ref{thm:relaxed_wellposedness}(iii). Bounded-cylinder comparison
\cite[Thm.~8.2]{CrandallIshiiLions1992} therefore yields
\(\underline U\le v^-_\eps\) and \(v^+_\eps\le\overline U\) on \(Q\);
the coefficients are bounded on \(\bar Q\) and the Hamiltonian is
Lipschitz in \(p\). The ordering \eqref{eq:value_bracketing} gives
the remaining value inequalities; Proposition~\ref{prop:value_level_reachability}
gives the set inclusions.
\end{proof}

A verified lateral shortfall \(c_+\) is corrected by adding \(c_+\) to the
upper candidate, without changing its HJI residual; the lower side is
analogous. Residual inequalities must hold throughout the cylinder, as
verified by outward-rounded interval arithmetic or a grid with a verified
remainder bound. An empirical loss or sampled maximum alone is insufficient.

\subsection{Boundary-compatible candidates}\label{subsec:boundary_candidates}
The lateral inequalities can be enforced by construction while retaining
flexibility in the interior. Let
\(\zeta,\beta\in C^\infty(\R^d;[0,1])\) with \(\zeta=1\), \(\beta=0\)
on a core region containing the states of interest, and \(\zeta=0\),
\(\beta=1\) on a neighborhood of \(\partial\Omega\).  With the hard
terminal ansatz \eqref{eq:hard_terminal_ansatz}, define
\[\begin{aligned}
    \overline u
    &=
    \widetilde g^+_\eps(x)
    +\frac{T-t}{T}\,\zeta(x)\,\mathcal N_{\theta^+}(t,x)
    +\beta(x),
    \\
    \underline u
    &=
    \widetilde g^-_\eps(x)
    +\frac{T-t}{T}\,\zeta(x)\,\mathcal N_{\theta^-}(t,x)
    -\beta(x).
\end{aligned}\]

{Assume
\(\widetilde g^-_\eps\le g^-_\eps\) and \(\widetilde g^+_\eps\ge g^+_\eps\) on \(\Omega\), and choose \(\Omega\) so that \(\widetilde g^\pm_\eps=0\) near \(\partial\Omega\). Then the terminal and lateral inequalities of Corollary~\ref{cor:localized_certificate} hold by construction, leaving only \(\eta_\pm\) to be verified. Alternatively, the lateral inequalities may be verified directly for other hard terminal ansatzes. On the compact cylinder \(\bar Q\), the residual suprema can be bounded by interval arithmetic over a finite cover or by a finite-grid evaluation with a verified Lipschitz remainder. Section~\ref{subsec:certified_1d} verifies the lateral inequalities and
residual bounds in one dimension by interval arithmetic with outward rounding on adaptive finite covers. The terms \(\pm\beta\) conservatively widen the bracket only outside the core region.}

\section{Neural candidate generation}\label{sec:PINNPI}
One-sided policy iteration motivates the candidate construction; neural
networks approximate the resulting policy-evaluation equations.

\subsection{Policy evaluation and improvement}
For disturbance and control feedbacks \(\beta\) and \(\alpha\), define
\[
 H_\beta(t,x,p)=\inf_a f(t,x,a,\beta(t,x))\cdot p,
 \qquad H^\alpha(t,x,p)=\sup_b f(t,x,\alpha(t,x),b)\cdot p.
\]
The Isaacs condition gives
\begin{equation}\label{eq:policy_hamiltonian_ordering}
 H_\beta\le H\le H^\alpha.
\end{equation}
Exact evaluation with \(H_\beta\) and terminal value \(g^-_\eps\)
yields a lower policy bound; evaluation with \(H^\alpha\) and terminal value
\(g^+_\eps\) yields an upper bound. Proposition~\ref{prop:policy_wise_bracketing}
in Appendix~\ref{app:howard} records the comparison argument.

\subsection{Exact policy iteration}
\label{subsec:exact_one_sided_howard}

The exact iterations use one-sided policy evaluation and greedy updates.

Choose an initial disturbance feedback policy \(\beta_0\).  Given
\(\beta_n\), solve the lower one-sided policy-evaluation equation
\begin{equation}
\label{eq:lower_howard_eval}
\begin{cases}
    \partial_t L_n
    +
    H_{\beta_n}(t,x,\nabla_x L_n)
    +
    \dfrac12\Tr(\Sigma D^2_{xx}L_n)
    =0,
    & (t,x)\in[0,T)\times\R^d,\\[1mm]
    L_n(T,x)=g^-_\eps(x),
    & x\in\R^d.
\end{cases}
\end{equation}
Then update the disturbance policy by the greedy selector
\begin{equation}
\label{eq:lower_howard_improve}
    \beta_{n+1}(t,x)
    \in
    \argmax_{b\in B}
    \inf_{a\in A}
    f(t,x,a,b)\cdot\nabla_x L_n(t,x).
\end{equation}
Equivalently,
\[
    H_{\beta_{n+1}}(t,x,\nabla_x L_n)
    =
    H(t,x,\nabla_x L_n).
\]

Choose an initial control feedback policy \(\alpha_0\).  Given
\(\alpha_n\), solve the upper one-sided policy-evaluation equation
\begin{equation}
\label{eq:upper_howard_eval}
\begin{cases}
    \partial_t U_n
    +
    H^{\alpha_n}(t,x,\nabla_x U_n)
    +
    \dfrac12\Tr(\Sigma D^2_{xx}U_n)
    =0,
    & (t,x)\in[0,T)\times\R^d,\\[1mm]
    U_n(T,x)=g^+_\eps(x),
    & x\in\R^d.
\end{cases}
\end{equation}
Then update the control policy by the greedy selector
\begin{equation}
\label{eq:upper_howard_improve}
    \alpha_{n+1}(t,x)
    \in
    \argmin_{a\in A}
    \sup_{b\in B}
    f(t,x,a,b)\cdot\nabla_x U_n(t,x).
\end{equation}
Equivalently,
\[
    H^{\alpha_{n+1}}(t,x,\nabla_x U_n)
    =
    H(t,x,\nabla_x U_n).
\]

\begin{proposition}
\label{prop:monotone_policy_certificates}
Assume bounded classical policy-evaluation solutions in
\eqref{eq:lower_howard_eval} and \eqref{eq:upper_howard_eval}, with
continuous terminal attainment and comparison against bounded classical
sub- and supersolutions. Assume also that the greedy selectors
\eqref{eq:lower_howard_improve} and \eqref{eq:upper_howard_improve}
exist.  Then
\[
    L_n
    \le
    L_{n+1}
    \le
    v^-_\eps
    \le
    v
    \le
    v^+_\eps
    \le
    U_{n+1}
    \le
    U_n
\]
on \([0,T]\times\R^d\).
\end{proposition}

\begin{proof}
The greedy lower update gives
\[
 H_{\beta_{n+1}}(t,x,\nabla L_n)=H(t,x,\nabla L_n)
 \ge H_{\beta_n}(t,x,\nabla L_n).
\]
Thus \(L_n\) is a subsolution of the next policy equation, with the same
terminal value, and comparison gives \(L_n\le L_{n+1}\). The identity
\(\mathcal H[L_n]=H-H_{\beta_n}\ge0\) gives \(L_n\le v^-_\eps\).
For the upper update, the reversed inequalities make \(U_n\) a
supersolution of the next equation, so \(U_{n+1}\le U_n\), while
\(\mathcal H[U_n]=H-H^{\alpha_n}\le0\) gives \(v^+_\eps\le U_n\).
Apply these bounds at both indices and use
\eqref{eq:value_bracketing}.
\end{proof}

\subsection{Convergence}
\label{subsec:exact_howard_convergence}

Identifying the monotone limits with the envelope values requires
compactness of the exact policy iterates and their derivatives.

\begin{assumption}
\label{ass:exact_howard_regularity}
For each fixed \(\eps>0\), the exact lower and upper one-sided
policy-evaluation equations
\eqref{eq:lower_howard_eval} and \eqref{eq:upper_howard_eval} admit
bounded classical solutions
\(L_n,U_n\in C([0,T]\times\R^d)\cap C^{1,2}([0,T)\times\R^d)\) that
attain their terminal data pointwise and continuously up to \(t=T\),
\(L_n(T,\cdot)=g^-_\eps\) and \(U_n(T,\cdot)=g^+_\eps\), and the
comparison principle against bounded classical sub- and supersolutions
with ordered terminal data holds for these policy-evaluation problems.
The greedy selectors in
\eqref{eq:lower_howard_improve} and \eqref{eq:upper_howard_improve}
exist as measurable feedback policies.

Moreover, for every compact \(K\subset\R^d\) and every
\(\tau\in(0,T]\), the exact policy-evaluation solutions satisfy uniform
local gradient bounds:
\[
    \sup_n
    \|\nabla_x L_n\|_{L^\infty([0,T-\tau]\times K)}
    +
    \sup_n
    \|\nabla_x U_n\|_{L^\infty([0,T-\tau]\times K)}
    <\infty .
\]
The sequences \(\{L_n\}_{n\ge0}\) and \(\{U_n\}_{n\ge0}\) are locally
precompact in \(C^{1,2}([0,T-\tau]\times K)\) for every compact
\(K\subset\R^d\) and every \(\tau\in(0,T]\).
\end{assumption}

Assumption~\ref{ass:exact_howard_regularity} concerns exact PDE solutions,
not network smoothness. The benchmark bang--bang selectors are discontinuous
at switching surfaces, so their required regularity and compactness need
separate estimates. Related controlled-diffusion convergence results are
given in \cite{KerimkulovSiskaSzpruch2020}.

\begin{theorem}
\label{thm:exact_howard_convergence}
Assume Assumption~\ref{ass:exact_howard_regularity}.  Then the exact
one-sided Howard iterates satisfy
\[
    L_n \uparrow v^-_\eps,
    \qquad
    U_n \downarrow v^+_\eps
\]
locally uniformly on \([0,T)\times\R^d\), with
\(L_n\le v\le U_n\) for every \(n\). Under
Assumption~\ref{ass:tube_estimate}, both envelope limits converge locally
uniformly to \(v\) as \(\eps\downarrow0\).
\end{theorem}
The proof is in Appendix~\ref{app:howard}.

\subsection{Policy error bounds}
\label{subsec:stagewise_policy_gaps}

For the exact iterates, suppose the following quantities are finite:
\[\begin{aligned}
 r_n^-&:=\sup_{(t,x)\in[0,T)\times\R^d}\mathcal H[L_n](t,x),\\
 r_n^+&:=\sup_{(t,x)\in[0,T)\times\R^d}
             \bigl(-\mathcal H[U_n](t,x)\bigr).
\end{aligned}\]
Both are nonnegative because \(H_{\beta_n}\le H\le H^{\alpha_n}\).

\begin{corollary}
\label{cor:stagewise_policy_certificate}
For every exact lower and upper Howard iterate,
\[\begin{aligned}
 0&\le v^-_\eps(t,x)-L_n(t,x)\le r_n^-(T-t),\\
 0&\le U_n(t,x)-v^+_\eps(t,x)\le r_n^+(T-t).
\end{aligned}\]
Consequently, \(L_n(t,x)\le v(t,x)\le U_n(t,x)\). Under
Assumption~\ref{ass:tube_estimate}, for \((t,x)\in[0,T-\tau]\times K\),
\[\begin{aligned}
 0&\le v(t,x)-L_n(t,x)\le r_n^-(T-t)+\omega_{K,\tau}(\eps),\\
 0&\le U_n(t,x)-v(t,x)\le r_n^+(T-t)+\omega_{K,\tau}(\eps).
\end{aligned}\]
\end{corollary}

The proof is in Appendix~\ref{app:howard}. The terms separate exact policy
error from terminal smoothing; bounded-domain certification additionally
requires lateral ordering.

\subsection{Neural policy evaluation}\label{subsec:pinn_policy_eval}
A hard terminal ansatz fixes \(u_\theta(T,x)=\psi(x)\); for example,
\begin{equation}\label{eq:hard_terminal_ansatz}
 u_\theta(t,x)=\psi(x)+\frac{T-t}{T}\mathcal N_\theta(t,x).
\end{equation}
The benchmark-specific parametrizations are stated in
Section~\ref{sec:impl_details}. For a frozen policy Hamiltonian
\(H_n\in\{H_{\beta_n},H^{\alpha_n}\}\), train with
\[
 \begin{gathered}
 \mathcal R_{\theta,n}=\partial_tu_\theta
 +H_n(t,x,\nabla_xu_\theta)+\tfrac12\Tr(\Sigma D^2_{xx}u_\theta),
 \\
 \mathcal J_n(\theta)=\frac1N\sum_{i=1}^N|\mathcal R_{\theta,n}(t_i,x_i)|^2.
 \end{gathered}
\]
Automatic differentiation supplies the derivatives. Alternatively, freezing
both policies gives a linear evaluation with
\(q^{\alpha_n,\beta_n}=f(t,x,\alpha_n,\beta_n)\cdot p\).
Even an exact solution of this linear equation has full HJI residual
\(H-q^{\alpha_n,\beta_n}\), whose sign is not prescribed. Thus the final
nonlinear residual and the terminal and lateral inequalities must still
be verified. The policy-gap identities are given in Appendix~\ref{app:howard}.

For certification, \(\psi\) bounds the indicator from the appropriate side;
the higher-dimensional experiments use a central smoothed payoff.
Training and sampling settings are given in Section~\ref{sec:impl_details}.
Write \(v_M(t,x;\theta_M)\) for the final candidate after \(M\) policy
updates, and \(v_R\) for the numerical reference used in each comparison.
\begin{algorithm}[t]
\caption{Candidate generation and certification}
\label{alg:pinn_howard}
\begin{algorithmic}[1]
\Require Terminal payoff \(\psi\), side \(s\in\{-,+,0\}\), initial policy, number of iterations \(M\)
\For{\(n=0,\dots,M-1\)}
    \State \textbf{Policy evaluation:} train \(u_{\theta_n}\) with
    \(u_{\theta_n}(T,\cdot)=\psi\) by minimizing the appropriate
    policy-evaluation residual.
    \State \textbf{Policy improvement:} update the active feedback
    policy by the greedy Hamiltonian selector using
    \(\nabla_x u_{\theta_n}\).
\EndFor
\State Set \(u_\theta\) to the final trained candidate.
\State For certification, verify terminal and lateral ordering and the one-sided bound for the
full HJI residual \(\mathcal H[u_\theta]\).
\State \Return candidate value function \(u_\theta\), and, when verified,
the corresponding residual margin.
\end{algorithmic}
\end{algorithm}

Here \(s=-,+\) denotes the one-sided envelope problems and \(s=0\) the
central-payoff approximation.

\section{Numerical evidence}\label{sec:experiments}
The Air3D and quadrotor benchmarks use
\(H_\eps(z)=\tfrac12(1+\tanh(z/\eps))\) to smooth target factors.
No finite-difference method (FDM) values enter the training loss.
Relative \(L^2\) errors are discrete Euclidean norm ratios on the stated
comparison grids. Additional plots and reference
details appear in Appendix~\ref{app:benchmarks}.

\subsection{Training and sampling}\label{sec:impl_details}
The archived Air3D and 10D implementations use
\[
 v_\theta(t,x)=(1-t)\mathcal N_\theta(t,x)+t g_\varepsilon(x),
 \qquad T=1.
\]
The reduced 4D implementation instead uses
\[
 v_\theta(t,x)=g_\varepsilon(x)+(T-t)\mathcal N_\theta(t,x),
 \qquad T=2.5.
\]
Both enforce the terminal condition. The scalar ansatz is given in
Section~\ref{subsec:certified_1d}. All networks use sine activations.

The higher-dimensional runs use Adam with learning rate \(10^{-3}\)
and pointwise minimax policy updates based on \(\nabla_xv\).
Let \(E\) be the maximum epochs per policy evaluation and \(N\) the
interior collocation count in Table~\ref{tab:training_config}.
Collocation sets are refreshed initially and every \(E/10\) epochs.
Air3D and 10D evaluations stop early if the relative loss change over
100 epochs falls below \(10^{-6}\), using denominator
\(\max(|\mathcal J_{\rm old}|,10^{-8})\).

Table~\ref{tab:training_config} gives the architectures and training budgets
for the higher-dimensional benchmarks. The 10D uniform-sampling baseline
uses the same architecture and nominal budget as the adaptive run.

\begin{table}[htbp]
\centering
\caption{PINN--PI training configurations; epoch counts are per-policy caps.}
\label{tab:training_config}
\footnotesize
\setlength{\tabcolsep}{3pt}
\renewcommand{\arraystretch}{1.0}
\begin{tabular}{lccc}
\hline
 & Air3D & 4D Quad & 10D Quad \\
\hline
Net. & 4--64$\times$4--1 & 5--64$\times$4--1 & 11--32$\times$4--1 \\
Params & 12,865 & 12,929 & 3,585 \\
Epoch cap/PI & 1k & 10k & 20k \\
PI iters & 1k & 500 & 500 \\
Colloc. & 4k & 50k & 50k \\
Boundary & 1k ($x_3$) & -- & -- \\
Sampling & uniform & adaptive & adaptive \\
\hline
\end{tabular}
\end{table}

Adaptive sampling concentrates points near the transition region. Draw \(M_c\) candidates
\[
    \{(\widetilde t^m,\widetilde x^m)\}_{m=1}^{M_c}
\]
uniformly in the computational space--time domain.  Using the current
network predictor, set \(\ell_m=\min\{1,\max\{0,
v_\theta(\widetilde t^m,\widetilde x^m)\}\}\). We draw \(N\) points
with replacement, with probability
\[
    \mathbb P\!\left((t,x)=(\widetilde t^m,\widetilde x^m)\right)
    \propto
    \max\{w(\ell_m),10^{-8}\}.
\]
The weight is
\[
    w(\ell)
    =
    (1-\gamma)\bigl(4\ell(1-\ell)\bigr)^\beta
    +
    \gamma\,\ell^p
    +
    \alpha\,(1-\ell)^p,
    \qquad
    \ell\in[0,1].
\]
For the 4D example, we set
\[
    (\gamma,\beta,p,\alpha)=(0.2,1.0,1.0,0.075).
\]
For the 10D example, the corresponding parameters are
\[
    (\gamma,\beta,p,\alpha)=(0.3,1.0,1.0,0.025),
\]
with \(g_\eps\) replacing \(v_\theta\) in the sampling weights during
the first 50 policy iterations.  In both examples, we draw \(M_c=8N\)
candidates before resampling. The first weight favors the transition region;
the remaining weights retain samples near \(v\approx1\) and \(v\approx0\).

\subsection{Scalar certification}
\label{subsec:certified_1d}

Consider the scalar game
\[
    dX_s=(a_s+b_s)\,ds+\sigma\,dW_s,
    \qquad
    a_s\in[-\bar a,\bar a],\quad b_s\in[-\bar b,\bar b],
\]
with \(\bar a=0.5\), \(\bar b=1.0\), \(\sigma=0.5\), horizon \(T=1\), target \(\mathcal G=[-0.5,0.5]\), and verification domain {\(\Omega=(-3.3,3.3)\)}. The Isaacs Hamiltonian is
\[
    H(p)
    =
    \sup_{|b|\le\bar b}\,\inf_{|a|\le\bar a}\,(a+b)p
    =
    (\bar b-\bar a)\,|p|.
\]
Both optimization orders coincide, satisfying Assumption~\ref{ass:isaacs}. With \(\rho(x)=0.5-|x|\), the one-sided envelopes \eqref{eq:one_sided_terminal_envelopes} with a smooth cutoff \(\chi\) and {\(\eps=0.4\)} are themselves smooth, since \(\chi(\rho/\eps)\) is locally constant near the kink of \(\rho\) at \(x=0\).
The implemented cutoff is
\(\chi(q)=[1+\exp(1/q-1/(1-q))]^{-1}\) for \(0<q<1\),
extended by 0 for \(q\le0\) and 1 for \(q\ge1\).

{We parameterize \(\underline u=u_{\theta^-}\) and \(\overline u=u_{\theta^+}\) by the hard terminal ansatz
\[
    u_{\theta^\pm}(t,x)
    =
    (1-t)\mathcal N_{\theta^\pm,\mathrm{even}}(t,x)
    +
    tG_\pm(x),
\]
where
\[
    \mathcal N_{\theta,\mathrm{even}}(t,x)
    :=
    \frac12\bigl(
    \mathcal N_\theta(t,x)+\mathcal N_\theta(t,-x)
    \bigr),
    \qquad
    G_-:=g^-_\eps-\beta,\quad
    G_+:=g^+_\eps+\beta .
\]
Here \( \beta(x)=\chi((|x|-2.5)/0.8) \), so that \(\beta=0\) on the core interval \([-2.5,2.5]\) and \(\beta=1\) on \(\partial\Omega\). The terminal inequalities therefore hold by construction, while the lateral inequalities are verified separately after training.

Each scalar candidate uses a \(2\)--\(32\times4\)--\(1\) sine network
with 3,297 parameters, trained for 500 policy evaluations of 16,000 steps
each. Adam uses cosine learning-rate decay from \(10^{-3}\) to
\(10^{-5}\) and gradient-norm clipping at 10. The residual batch contains
17,068 interior points and two 1,024-point terminal batches; 512 boundary
time samples supply the lateral penalty. Half of each random sample block
is replaced every 1,000 steps, retaining the deterministic terminal batch.
The archived configurations specify the
terminal sampling and random seeds.

At policy iteration \(n\), the frozen policy
\[
    q_n^\pm
    =
    \operatorname{sign}\bigl(\partial_xu_{\theta_n^\pm}\bigr)
\]
gives the policy-evaluation residual
\[
    \mathcal R_n^\pm[u]
    =
    \partial_tu
    +0.5q_n^\pm\partial_xu
    +0.125\partial_{xx}u .
\]
The loss adds to the empirical mean square of \(\mathcal R_n^\pm\)
the mean squared largest 1\% of wrong-sign violations with weight 0.1
and a lateral boundary penalty with weight 100.
Interval arithmetic with outward rounding verifies the lower lateral
inequality without correction, corresponding to \(c_-=0\), and bounds the lateral shortfall of the upper candidate by less than \(10^{-3}\). We therefore define
\(
\overline u_{\rm corr}:=\overline u+c_+,
\ c_+=10^{-3}
\).
Since the HJI operator is invariant under constant shifts, this correction preserves both the upper terminal inequality and the full HJI residual.}

{The full HJI residuals \(r^-:=\mathcal H[\underline u]\) and \(r^+:=\mathcal H[\overline u_{\rm corr}]\) are bounded using outward-rounded interval arithmetic over an adaptive finite cover of \(\bar Q=[0,T]\times[-3.3,3.3]\). By symmetry, it suffices to consider \(x\ge0\). Writing \(z_+=\max\{z,0\}\) and \(z_-=\max\{-z,0\}\), the interval bounds give
\[
    \bigl(r^-\bigr)_-<0.1,
    \qquad
    \bigl(r^+\bigr)_+<0.1
    \qquad\text{on }\bar Q,
\]
so we take \(\eta_-=\eta_+=0.1\).}

The accepted covers contain 202,081 lower-candidate and 239,406
upper-candidate cells on \([0,1]\times[0,3.3]\). The public audit checks
each stored residual enclosure against the cap and verifies gap-free,
nonoverlapping coverage using exact rational endpoints. Reflection covers
the negative half. Fresh enclosure recomputation is documented separately
in Appendix~\ref{app:computational_audit}.

Corollary~\ref{cor:localized_certificate} gives
\(\underline u-\eta_-(T-t)\le v\le\overline u_{\rm corr}+\eta_+(T-t)\)
on \(Q\). {At \(t=0\), \(\delta=0.5\), and \(\xi=0.01\), interval
classification yields
\[
\begin{aligned}
[-0.546,0.546]
&\subset R_{0.5}(0)\cap\Omega \\
&\subset (-3.3,-2.258]\cup[-1.332,1.332] \\
&\qquad\cup[2.258,3.3).
\end{aligned}
\]
The outer boundary components arise from the lateral upper construction.
Independent monotone FDM references satisfy the bracket at the \(t=0\)
grid nodes. They use terminal data \(G_\pm\), Dirichlet values \(-1\)
and \(1\), respectively, and spacing \(0.00125\). These bounded-domain
references include the terminal and lateral modifications.}

At \(t=0\), the widest bracket on the core grid \([-2.5,2.5]\) occurs
at \(x=2.5\). Its width \(0.83428\) decomposes into a bounded-domain
reference gap of \(0.63327\), a signed candidate--reference contribution of
\(0.00101\), and residual correction \(0.20000\). The reference gap includes
terminal smoothing and lateral modification, so improved candidate--FDM
agreement alone would not substantially narrow this bracket. These values
are numerical diagnostics, not interval bounds on the whole-space
envelope gap.

\begin{figure}[htbp]
    \centering
    \includegraphics[width=\columnwidth]{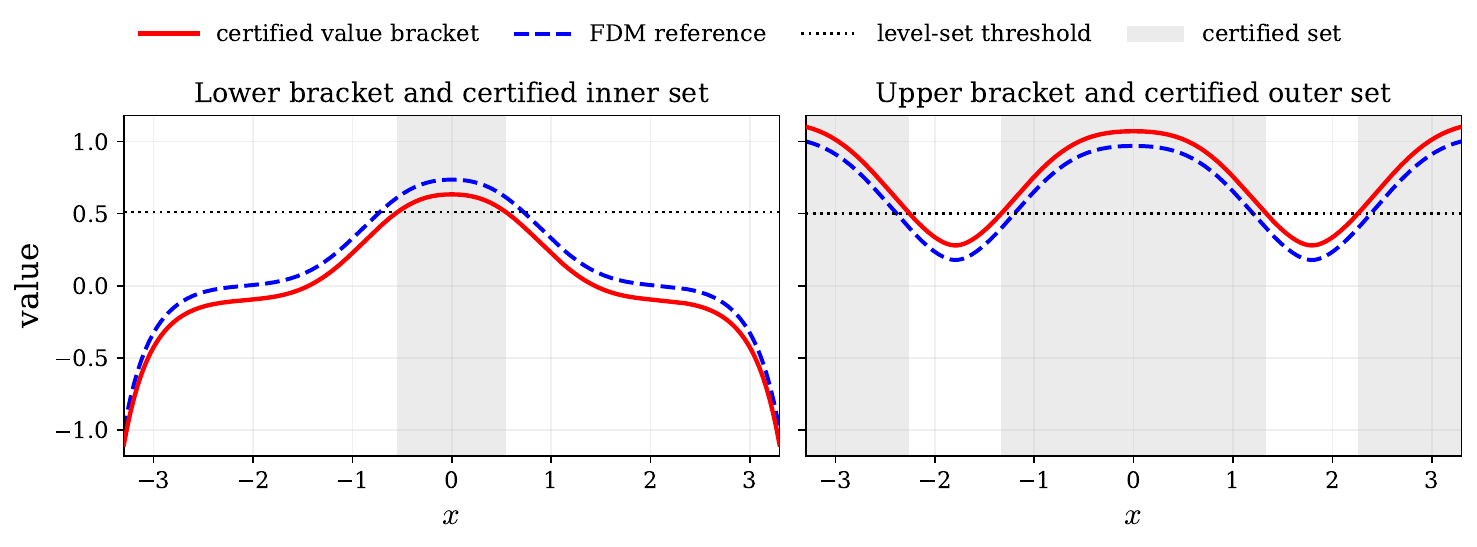}
    \caption{\textrm{Certified comparison at \(t=0\) for \(\delta=0.5\), with \(\xi=0.01\) for the inner approximation. The red curves are the certified lower and upper bounds with \(\eta_-=\eta_+=0.1\), the blue dashed curves are the independent monotone FDM references, and the shaded regions are the certified inner and outer approximations.  The boundary components of the outer approximation result from the enforced lateral upper bound.}} \label{fig:1d_certified_bracket}
\end{figure}

\subsection{Air3D}\label{subsec:air3d}
\begin{figure}[htbp]
\centering
\includegraphics[width=\columnwidth]{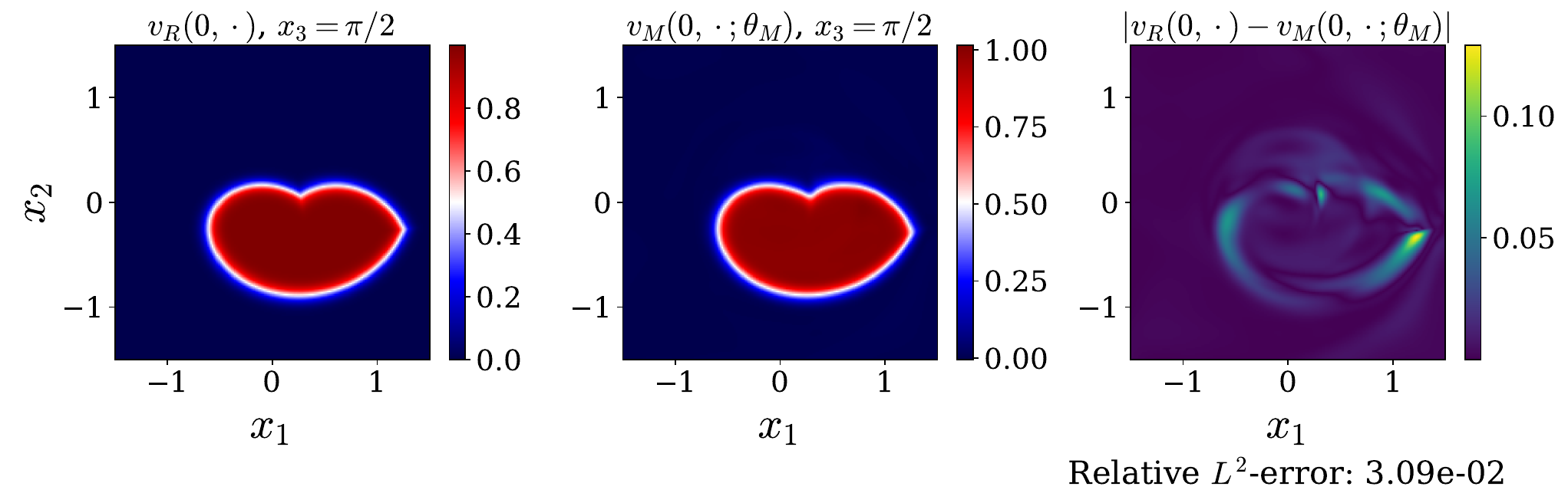}
\caption{Air3D example.
Left: FDM reference solution for the smoothed terminal payoff.
Middle: PINN--PI solution. Right: absolute error.}
\label{fig:air3d_pinn_fdm_error}
\end{figure}

The stochastic Air3D collision-avoidance benchmark
\cite{bansal2021deepreach} has state \(x=(x_1,x_2,x_3)\), where the
relative heading \(x_3\) is \(2\pi\)-periodic. Its drift is
\[
 f(x,\omega_e,\omega_p)=
 \begin{pmatrix}
 -v_e+v_p\cos x_3+\omega_e x_2\\
 v_p\sin x_3-\omega_e x_1\\
 \omega_p-\omega_e
 \end{pmatrix},\qquad
 \sigma_\kappa=\kappa
 \begin{pmatrix}
 0.05&0.02&0.01\\0.02&0.03&0\\0.01&0&0.10
 \end{pmatrix}.
\]
We use \(v_e=v_p=0.5\), \(|\omega_e|,|\omega_p|\le\omega_{\max}=2\),
\(T=1\), and computational domain \([-1.5,1.5]^2\times[-\pi,\pi]\).
For gradient variable \(p\), the full Isaacs Hamiltonian is
\[
 \begin{aligned}
 H(x,p)&=(-v_e+v_p\cos x_3)p_1+v_p\sin x_3\,p_2\\
 &\quad+\omega_{\max}|p_3|
 -\omega_{\max}|x_2p_1-x_1p_2-p_3|.
 \end{aligned}
\]
The target is the position disk of radius 0.75, and
\(g_\eps(x)=H_\eps(0.75-\sqrt{x_1^2+x_2^2})\) with \(\eps=0.05\).
Figure~\ref{fig:air3d_pinn_fdm_error} compares the learned candidate with
FDM at \(t=0\), \(x_3=\pi/2\), and unit diffusion scale. The relative
\(L^2\) discrepancy on this slice is approximately \(3.09\times10^{-2}\).

We also evaluate terminal hitting probabilities under the learned feedback
pair for diffusion scales \(\kappa=0.5,1,5\). For each scale, 1,000 initial
states (seeds 0--999) are selected from the learned value band
\(v_\kappa(0,x)\in[0.4,0.5]\) on a \(301\times301\) position grid in
\([-1.5,1.5]^2\), with \(x_3=\pi/2\). Each state is evaluated using
1,000 Euler--Maruyama paths with \(\Delta t=0.01\). The estimate is
\(\widehat p_\kappa(x)=\#\{X_T\in\mathcal G\}/1000\).
Table~\ref{tab:air3d_mc_kappa} reports its distribution across states.
Feedback formulas and representative trajectories appear in
Appendix~\ref{app:benchmarks}.
\begin{table}[htbp]
    \centering
    \caption{Air3D Monte Carlo diagnostic under diffusion
    scaling.  For each \(\kappa\), we sample \(1000\) initial states at
    \(t=0\) from the learned value band
    \(v_\kappa(0,x)\in[0.4,0.5]\).  For each state \(x\), we estimate
    \(\widehat p_\kappa(x)=\#\{X_T\in\mathcal G\}/1000\) using \(1000\)
    Euler--Maruyama rollouts. The table summarizes terminal-event
    probabilities under the learned feedback pair.}
    \label{tab:air3d_mc_kappa}
    \begin{tabular}{ccccc}
        \hline
        \(\kappa\) & mean (std) & min & median & max \\
        \hline
        0.5 & 0.1507 (8.51e-02) & 0.0150 & 0.1360 & 0.4580 \\
        1.0 & 0.2974 (8.84e-02) & 0.0970 & 0.3130 & 0.5360 \\
        5.0 & 0.4144 (4.53e-02) & 0.2850 & 0.4150 & 0.5480 \\
        \hline
    \end{tabular}
\end{table}

The means 0.1507 and 0.2974 at the smaller diffusion scales lie below the
selected band; agreement is closer at \(\kappa=5\). The network approximates
a smoothed minimax value, whereas rollouts evaluate an indicator under one
learned feedback pair. Their difference includes smoothing, feedback,
PDE, and time-discretization effects, which this experiment does not isolate.

\subsection{Reduced 4D quadrotor}\label{subsec:quad4d}
\begin{figure}[htbp]
\centering
\includegraphics[width=\columnwidth]{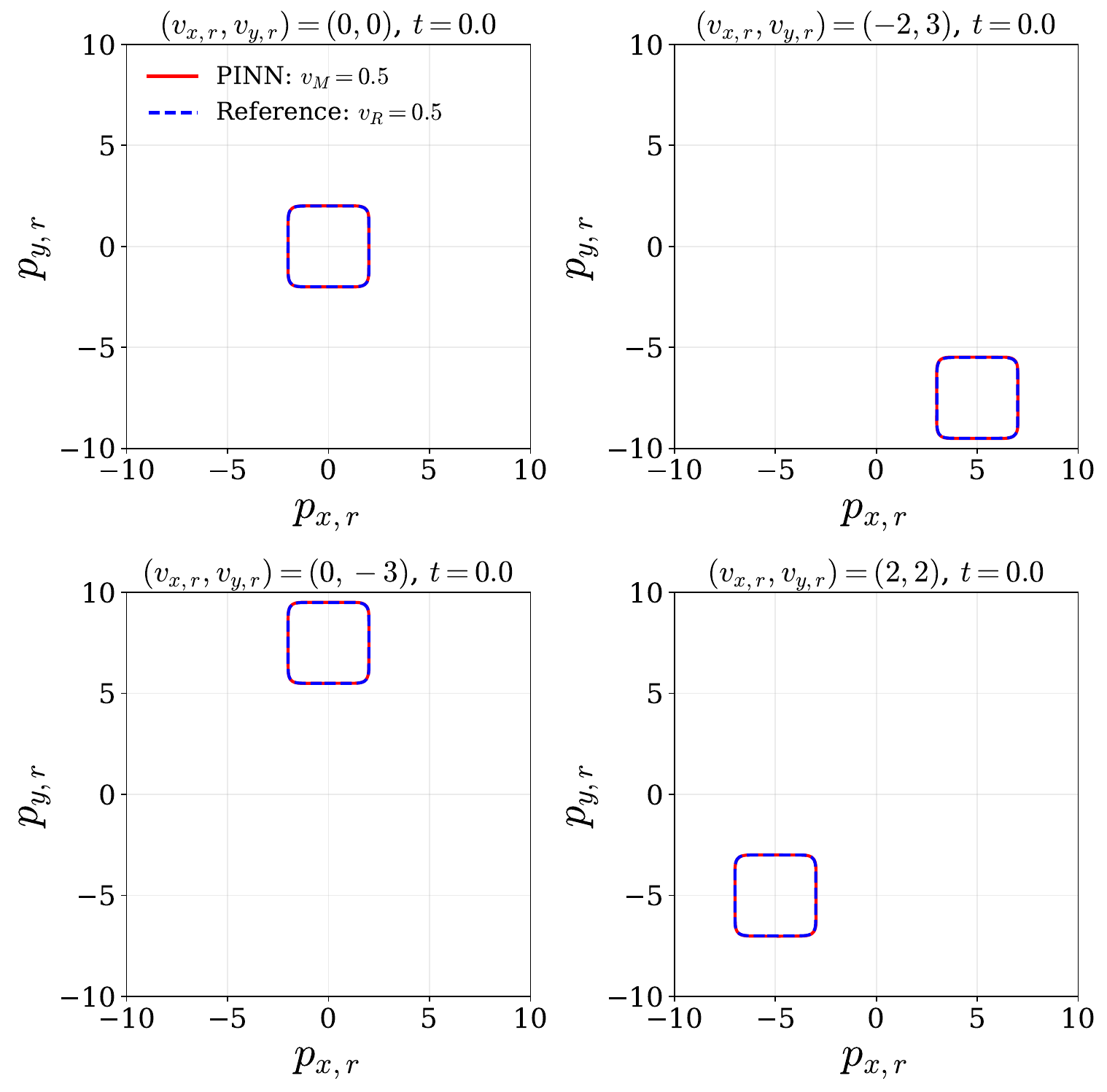}
\caption{Reduced 4D quadrotor: value-function slices at \(t=0\).
Each panel shows the \((p_{x,r},p_{y,r})\)-slice at a fixed
\((v_{x,r},v_{y,r})\).  Red solid: PINN--PI contour
\(v_M(t,x;\theta_M)=0.5\).  Blue dashed: processed finite-difference reference
contour \(v_R(t,x)=0.5\), constructed from the product-structured
subsystem reference in Appendix~\ref{app:decomp}.}
\label{fig:quad4d:t0:level05:slices}
\end{figure}

The reduced relative-position model \cite{chen2015exact} omits the
evader velocity states and the original state constraint. It has two
position--velocity blocks \(x^{(k)}=(p_r,v_r)\), each with dynamics
\[
 d\begin{pmatrix}p_r\\v_r\end{pmatrix}
 =\begin{pmatrix}v_r\\u-d\end{pmatrix}ds+\sigma_k\,dW_s^{(k)},
 \qquad \sigma_k=\begin{pmatrix}0.05&0.02\\0.02&0.07\end{pmatrix},
\]
where the Brownian blocks are independent, \(|u|\le u_{\max}\), and
\(|d|\le d_{\max}\). The computational domain is
\([-10,10]\times[-6,6]\times[-10,10]\times[-6,6]\). We use
\(T=2.5\), target half-widths \(d_x=d_y=2\), and terminal payoff
\[
 g_\eps(x)=H_\eps(d_x-|p_{x,r}|)H_\eps(d_y-|p_{y,r}|),\qquad\eps=0.2.
\]
The subsystem Hamiltonian is
\[
 H_k(x^{(k)},p)=v_rp_1+(d_{\max}-u_{\max})|p_2|.
\]
In the implemented case \(u_{\max}=d_{\max}=1\), the optimized inputs
cancel, giving a linear advection--diffusion benchmark.
We compare the
\((p_{x,r},p_{y,r})\) slices at velocities
\((0,0),(-2,3),(0,-3),(2,2)\); their relative \(L^2\) errors against the
processed FDM reference range from about 1.25\% to 2.03\%.
Figure~\ref{fig:quad4d:t0:level05:slices} shows the level-0.5 contours.
The processing retains a principal connected component before subsampling;
its exact definition and the absolute-error panels are in
Appendix~\ref{app:benchmarks}.

The cancellation also permits an independent Gaussian terminal-expectation
check. For \(h=T-t\) and \(Q=\sigma_k\sigma_k^\top\), each terminal
position is Gaussian with
\[
 \E[P_T]=p+hv,\qquad
 \operatorname{Var}(P_T)=hQ_{11}+h^2Q_{12}+\frac{h^3}{3}Q_{22}.
\]
Independence of the two blocks makes the comparator the product of their
Gaussian expectations of the smoothed target factors. Increasing the
quadrature order from 128 to 256 changes it by less than
\(1.8\times10^{-11}\) in maximum norm.
On the four \(201\times201\) slices, the unmodified FDM reference
differs from this comparator by relative \(L^2\) error below
\(8\times10^{-5}\); neural discrepancies range from
\(1.25\times10^{-2}\) to \(2.02\times10^{-2}\).
Quadrature and post-processing details are in Appendix~\ref{app:computational_audit}.

\subsection{10D quadrotor}\label{subsec:quad10d}
The stochastic 10D near-hover quadrotor benchmark
\cite{chen2016decomposition} has state
\[
    x=(p_x,v_x,\theta_x,\omega_x,
       p_y,v_y,\theta_y,\omega_y,
       p_z,v_z).
\]
The target set is
\[
    \mathcal G
    =
    \{|p_x|\le1,\ |p_y|\le1,\ |p_z|\le2.5\},
\]
and the smoothed terminal payoff is
\[
    g_\eps(x)
    =
    H_\eps(1-|p_x|)
    H_\eps(1-|p_y|)
    H_\eps(2.5-|p_z|),
    \qquad
    \eps=0.2.
\]

We add a block-diagonal diffusion term consistent with the decoupled
structure of the near-hover benchmark:
\[
    dX_s=f(X_s,a_s,b_s)\,ds+\sigma\,dW_s,
    \qquad
    \sigma=\diag(\sigma_x,\sigma_y,\sigma_z).
\]
The control and disturbance inputs are
\[
    a_s=(S_x,S_y,T_z),
    \qquad
    b_s=(d_x,d_y,d_z).
\]
The drift is the near-hover model of
\cite{chen2016decomposition}: with gravitational acceleration
\(g_{\rm grav}\) and positive model parameters \(d_0\), \(d_1\),
\(n_0\), \(k_T\), the deterministic part of the dynamics is
\[
\begin{aligned}
\dot p_x&=v_x+d_x,&\dot v_x&=g_{\rm grav}\tan\theta_x,\\
\dot\theta_x&=-d_1\theta_x+\omega_x,&\dot\omega_x&=-d_0\theta_x+n_0S_x,\\
\dot p_y&=v_y+d_y,&\dot v_y&=g_{\rm grav}\tan\theta_y,\\
\dot\theta_y&=-d_1\theta_y+\omega_y,&\dot\omega_y&=-d_0\theta_y+n_0S_y,\\
\dot p_z&=v_z+d_z,&\dot v_z&=k_TT_z-g_{\rm grav}.
\end{aligned}
\]
where \(S_x\) and \(S_y\) are the pitch and roll commands, \(T_z\) is
the vertical thrust, and \(d_x\), \(d_y\), \(d_z\) are the
disturbances.  On the computational domain,
\(|\theta_x|,|\theta_y|\le0.2\), so \(\tan\theta\) is smooth and
Lipschitz there; the global Lipschitz requirement of
Assumption~\ref{ass:dynamics} can be met by a globally Lipschitz extension outside
this range. Such an extension leaves the evaluated local operator unchanged;
it does not establish equality of global game values.

The implementation uses \(g_{\rm grav}=9.81\), \(d_0=n_0=10\),
\(d_1=8\), \(k_T=0.91\), \(S_{x,\max}=S_{y,\max}=10\pi/180\),
\(T_{z,\min}=0\), \(T_{z,\max}=2(9.81)\), and
\((\bar d_x,\bar d_y,\bar d_z)=(0.5,0.5,1.0)\).

The horizon is \(T=1\). Appendix~\ref{app:benchmarks} specifies the
computational box, diffusion matrices, and subsystem Hamiltonians.
The neural candidate uses all ten states.
Quantitative evaluation is on
the \((p_x,p_y)\) plane with \(p_z=0\),
\((v_x,v_y,v_z)=(-1.5,-1.8,1.2)\), and zero angles and angular velocities.
The FDM comparator is a subsystem product. Appendix~\ref{app:decomp}
establishes factorization for a restricted product game: each block uses
only its own Brownian filtration, with independent noise and no shared
randomization. Under these restrictions the reference has the form
\(v_R=V_1V_2V_3\) for the two horizontal four-state blocks and the
vertical two-state block. The factorization does not establish equality
with the unrestricted 10D game value; the reported errors retain this
qualification.

\begin{figure}[htbp]
    \centering
    \includegraphics[width=\columnwidth]{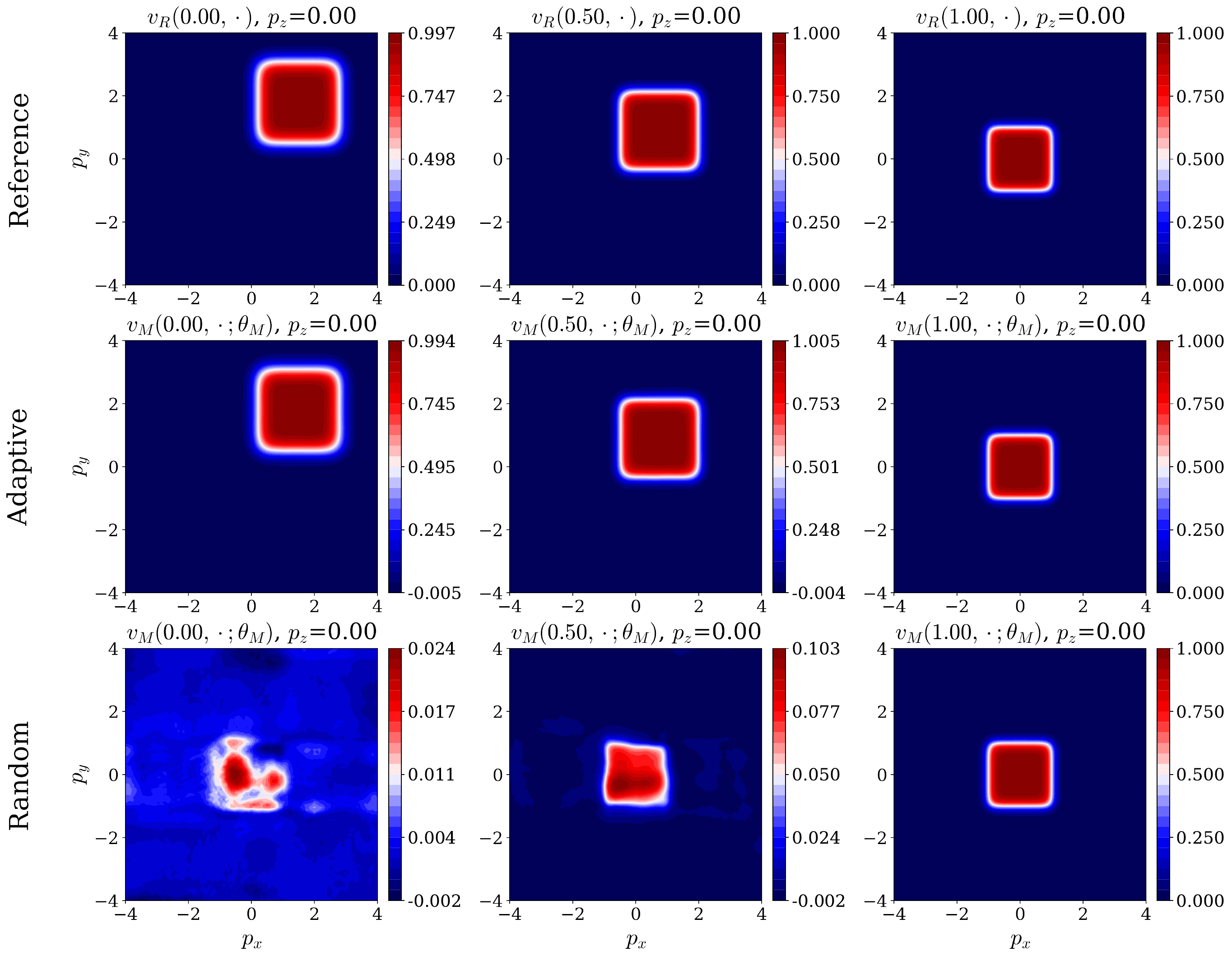}
    \caption{10D near-hover quadrotor: comparison of value-function
    slices under adaptive and random collocation sampling.  Columns
    correspond to \(t\in\{0,0.5,1.0\}\).  Rows show the FDM reference,
    the PINN--PI result with adaptive sampling, and the PINN--PI result
    with random sampling.  Each panel shows the \((p_x,p_y)\)-slice at
    \(p_z=0\) with the remaining states fixed.}
    \label{fig:nh10d_adaptive_random_t0_t05_t10}
\end{figure}

With the same network and nominal training budget, adaptive and uniform
collocation give relative errors of 1.15\% and 99.6\%, respectively, on
the evaluated slices at \(t=0\)
(Figure~\ref{fig:nh10d_adaptive_random_t0_t05_t10} and
Table~\ref{tab:nh10d_l2_time_adaptive_random}). This run supports
concentrating collocation points near the transition region.

\begin{table}[H]
    \centering
    \caption{10D near-hover quadrotor: relative \(L^2\) errors against
    the FDM reference at different times.  Errors are computed on the
    \((p_x,p_y)\)-slice at \(p_z=0\) with the remaining states fixed.
    These are slice diagnostics against the restricted-product reference.}
    \label{tab:nh10d_l2_time_adaptive_random}
    \setlength{\tabcolsep}{3.5pt}
    \renewcommand{\arraystretch}{1.0}
    \footnotesize
    \begin{tabular}{cccccc}
        \hline
        \(t\) & 0.00 & 0.25 & 0.50 & 0.75 & 1.00 \\
        \hline
        Adap. & 1.15e-02 & 1.01e-02 & 9.13e-03 & 8.69e-03 & 1.09e-07 \\
        Rand. & 9.96e-01 & 9.90e-01 & 9.72e-01 & 9.31e-01 & 1.09e-07 \\
        \hline
    \end{tabular}
\end{table}

\begin{samepage}
The terminal row measures reference and rounding consistency: the network
contribution vanishes at \(t=T\). Appendix~\ref{app:computational_audit}
records regeneration differences.\par
\end{samepage}

\FloatBarrier
\section{Conclusion}
The residual certificate turns smooth candidates into bounds on
strategy-defined reachability, as demonstrated by the scalar interval
verification and its nonempty inner set. The higher-dimensional experiments
support neural candidate generation, with adaptive collocation reaching
about one-percent relative error against the restricted-product reference
on the evaluated 10D slices.

\appendix
\renewcommand{\thesection}{\Alph{section}}
\numberwithin{equation}{section}
\numberwithin{figure}{section}
\numberwithin{table}{section}

\section{Comparison theory and discontinuous terminal data}\label{sec:HJI}

Payoff ordering and viscosity comparison justify the passage from
continuous terminal envelopes to the discontinuous target indicator.

\subsection{Dynamic programming and payoff ordering}
For bounded uniformly continuous \(\psi\), the standard stochastic-game
dynamic programming principle \cite{FlemingSouganidis1989,BuckdahnLi2008}
gives, at deterministic \(s\in[t,T]\),
\begin{equation}\label{eq:DPP_general_payoff}
 V^\psi(t,x)=\sup_{\gamma^0\in\Gamma_{t,s}}\inf_{a^0\in\mathcal A_{t,s}}
 \E[V^\psi(s,X_s^{t,x;a^0,\gamma^0[a^0]})].
\end{equation}
It is used only for continuous terminal payoffs. Indicator-payoff bounds
follow instead from monotonicity and the terminal-envelope construction.
\begin{lemma}
\label{lem:value_monotonicity}
Let \(\psi_1,\psi_2\in B_b(\R^d)\).  If
\[
    \psi_1(x)\le \psi_2(x),
    \qquad x\in\R^d,
\]
then
\[
    V^{\psi_1}(t,x)\le V^{\psi_2}(t,x),
    \qquad
    (t,x)\in[0,T]\times\R^d.
\]
Moreover,
\[
    \inf_{\R^d}\psi
    \le
    V^\psi(t,x)
    \le
    \sup_{\R^d}\psi.
\]
If \(c\in\R\), then
\[
    V^{\psi+c}(t,x)=V^\psi(t,x)+c.
\]
\end{lemma}

\begin{proof}
Expectation, infimum, and supremum preserve order and commute with addition
of a constant. Applying this to constant lower and upper payoffs gives
the bounds.
\end{proof}

For the indicator payoff the formal problem is
\begin{equation}\label{eq:HJI}
 \mathcal H[v]=0\quad(t<T),\qquad v(T,\cdot)=\1_{\mathcal G},
\end{equation}
with the terminal condition interpreted through semicontinuous traces.
\subsection{Relaxed viscosity formulation}

For a locally bounded function \(u\) on \([0,T)\times\R^d\), define its
upper and lower semicontinuous terminal envelopes by
\[
    u^*(T,x)
    :=
    \limsup_{\substack{(s,y)\to(T,x)\\ s<T}} u(s,y),
    \qquad
    u_*(T,x)
    :=
    \liminf_{\substack{(s,y)\to(T,x)\\ s<T}} u(s,y).
\]
Inside the cylinder \([0,T)\times\R^d\), \(u^*\) and \(u_*\) are defined
in the usual way.

For the indicator payoff \(g=\1_{\mathcal G}\), define
\[
    g^*(x):=\limsup_{y\to x}\1_{\mathcal G}(y),
    \qquad
    g_*(x):=\liminf_{y\to x}\1_{\mathcal G}(y).
\]
Since \(\mathcal G\) is closed,
\begin{equation}
\label{eq:indicator_envelopes}
    g^*(x)=\1_{\mathcal G}(x),
    \qquad
    g_*(x)=\1_{\operatorname{int}\mathcal G}(x).
\end{equation}
In particular, on \(\partial\mathcal G\),
\[
    g_*(x)=0,
    \qquad
    g^*(x)=1.
\]
The assigned terminal value is still \(v(T,x)=\1_{\mathcal G}(x)\).
At boundary points, it need not equal the limit from \(t<T\).

We use the following sign convention.  Since the equation is written as
\(\mathcal L u=0\), we apply the standard viscosity convention to the
proper operator \(-\mathcal L\).  Equivalently, a subsolution is tested
by the inequality \(\mathcal L\phi\ge0\) at local maxima of
\(u^*-\phi\), while a supersolution is tested by
\(\mathcal L\phi\le0\) at local minima of \(u_*-\phi\).

\begin{definition}
\label{def:relaxed_viscosity_solution}
A locally bounded function \(u\) on \([0,T)\times\R^d\) is a relaxed
viscosity subsolution of \eqref{eq:HJI} if \(u^*\) is a viscosity
subsolution of \(\mathcal L u=0\) on \([0,T)\times\R^d\) and satisfies
\[
    u^*(T,x)\le g^*(x),
    \qquad x\in\R^d.
\]
It is a relaxed viscosity supersolution if \(u_*\) is a viscosity
supersolution of \(\mathcal L u=0\) on \([0,T)\times\R^d\) and satisfies
\[
    u_*(T,x)\ge g_*(x),
    \qquad x\in\R^d.
\]
A relaxed viscosity solution is a locally bounded function that is both
a relaxed viscosity subsolution and a relaxed viscosity supersolution.
\end{definition}

The relaxed terminal condition is sharp for discontinuous terminal
payoffs.  At points \(x\in\operatorname{int}\mathcal G\), one expects
\(v(t,x)\to1\) as \(t\uparrow T\), while at points
\(x\notin\mathcal G\), one expects \(v(t,x)\to0\).  At boundary points
\(x\in\partial\mathcal G\), however, the limit need not equal
\(\1_{\mathcal G}(x)\).  For instance, in the presence of a
nondegenerate diffusion, the short-time probability of being in
\(\mathcal G\) from a boundary point may converge to a value strictly
between \(0\) and \(1\).

The following elementary estimate on the state process is the only
probabilistic ingredient needed to handle the terminal layer.  It
requires no dynamic programming and holds uniformly over all admissible
controls and disturbance processes.

\begin{lemma}
\label{lem:short_time_moment}
Assume Assumption~\ref{ass:dynamics}.  For every \(R>0\) there exists
\(C_R>0\) such that, for all \(t\in[0,T]\), all \(|x|\le R\), all
\(a\in\mathcal A_t\), and all \(b\in\mathcal B_t\),
\[
    \E\!\left[
        \sup_{s\in[t,T]}
        \bigl|X_s^{t,x;a,b}-x\bigr|^2
    \right]
    \le
    C_R\,(T-t).
\]
The same bound holds with \(b=\gamma[a]\) for any strategy
\(\gamma\in\Gamma_t\), since \(\gamma[a]\in\mathcal B_t\).
\end{lemma}

\begin{proof}
Write \(X=X^{t,x;a,b}\).  By the linear growth of \(f\) and \(\sigma\),
uniform over the compact sets \(A\times B\), the standard Gronwall
estimate gives \(\E[\sup_{s\in[t,T]}|X_s|^2]\le C(1+|x|^2)\), with \(C\)
depending only on the coefficient bounds and \(T\).  From
\(X_s-x=\int_t^s f(r,X_r,a_r,b_r)\,dr+\int_t^s\sigma(r,X_r)\,dW_r\), the
Cauchy--Schwarz inequality and Doob's maximal inequality yield
\[
    \begin{aligned}\E\Bigl[\sup_{s\in[t,T]}|X_s-x|^2\Bigr]
    &\le
    2(T-t)\int_t^T\E\bigl[|f(r,X_r,a_r,b_r)|^2\bigr]\,dr
    \\
    &\quad+
    8\int_t^T\E\bigl[|\sigma(r,X_r)|^2\bigr]\,dr,
\end{aligned}\]
and both integrands are bounded by \(C(1+\E[|X_r|^2])\le C'(1+R^2)\).
Hence the left-hand side is bounded by
\(C_R((T-t)^2+(T-t))\le C_R'(T-t)\).  All constants are independent of
\((a,b)\).
\end{proof}

\begin{theorem}
\label{thm:relaxed_wellposedness}
Assume Assumptions~\ref{ass:dynamics} and~\ref{ass:isaacs}.
\begin{enumerate}
\item[(i)] \emph{(Relaxed terminal inequalities for \(v\).)}  The game
value \(v\) defined in \eqref{eq:indicator_game_value} satisfies
\[
    v^*(T,x)\le g^*(x),
    \qquad
    v_*(T,x)\ge g_*(x),
    \qquad x\in\R^d.
\]
\item[(ii)] \emph{(Comparison with ordered terminal traces.)}  Let \(u\)
and \(w\) be bounded functions on \([0,T)\times\R^d\) such that \(u^*\)
is a viscosity subsolution and \(w_*\) is a viscosity supersolution of
\(\mathcal Lu=0\) on \([0,T)\times\R^d\), and suppose
\[
    u^*(T,x)\le w_*(T,x),
    \qquad x\in\R^d.
\]
Then \(u\le w\) on \([0,T)\times\R^d\).
\item[(iii)] \emph{(Continuous terminal data.)}  For every bounded
uniformly continuous \(\psi:\R^d\to\R\), the game value \(V^\psi\) is
continuous on \([0,T]\times\R^d\), attains its terminal datum locally
uniformly,
\[
    \sup_{|x|\le R}\bigl|V^\psi(t,x)-\psi(x)\bigr|
    \xrightarrow[t\uparrow T]{}0
    \qquad\text{for every }R>0,
\]
and is the unique bounded viscosity solution of
\eqref{eq:HJI_general_payoff} on \([0,T)\times\R^d\) with terminal datum
\(\psi\).  In particular, the envelope problems with terminal data
\(g^-_\eps\) and \(g^+_\eps\) defined in Section~\ref{sec:formulation} are well posed, and their
solutions are the game values \(v^\pm_\eps=V^{g^\pm_\eps}\).
\end{enumerate}
\end{theorem}

\begin{proof}
(i) Both inequalities follow from Lemma~\ref{lem:short_time_moment} and
Chebyshev's inequality; no dynamic programming principle for the
discontinuous payoff is needed.

For the first inequality, only points \(x\notin\mathcal G\) require an
argument, since \(g^*\equiv1\) on \(\mathcal G\) and \(v\le1\).  Fix
\(x\notin\mathcal G\) and set
\(\eta:=\min\{1,\tfrac12\operatorname{dist}(x,\mathcal G)\}>0\) (recall that
\(\mathcal G\) is closed).  For \(|y-x|\le\eta/2\) and \(t<T\), every
admissible pair \((a,\gamma)\) satisfies
\(\{X_T^{t,y;a,\gamma[a]}\in\mathcal G\}
\subset\{|X_T^{t,y;a,\gamma[a]}-y|\ge\eta\}\), whence, by Chebyshev's
inequality and Lemma~\ref{lem:short_time_moment} with \(R=|x|+1\),
\[
    \Prob\bigl(X_T^{t,y;a,\gamma[a]}\in\mathcal G\bigr)
    \le
    \frac{C_R\,(T-t)}{\eta^2}.
\]
The bound is uniform in \((a,\gamma)\), so
\(v(t,y)\le C_R(T-t)/\eta^2\), and letting \((t,y)\to(T,x)\) with
\(t<T\) gives \(v^*(T,x)\le0=g^*(x)\).

For the second inequality, only points
\(x\in\operatorname{int}\mathcal G\) require an argument.  Fix such an
\(x\) and set \(\eta:=\min\{1,\tfrac12\operatorname{dist}(x,\R^d\setminus
\operatorname{int}\mathcal G)\}>0\), with \(R=|x|+1\), so that
\(\bar B(x,2\eta)\subset\mathcal G\).  Fix any \(\bar b\in B\) and let
\(\bar\gamma\in\Gamma_t\) be the constant strategy
\(\bar\gamma[a]\equiv\bar b\).  For \(|y-x|\le\eta/2\), the inclusion
\(\{|X_T-y|\le\eta\}\subset\{X_T\in\mathcal G\}\) holds (if
\(|z-y|\le\eta\), then \(|z-x|\le\tfrac32\eta<2\eta\)), so for every
\(a\in\mathcal A_t\),
\[
    \Prob\bigl(X_T^{t,y;a,\bar\gamma[a]}\in\mathcal G\bigr)
    \ge
    1-\frac{C_R\,(T-t)}{\eta^2}.
\]
Hence
\(v(t,y)\ge\inf_{a}\Prob(X_T^{t,y;a,\bar\gamma[a]}\in\mathcal G)
\ge1-C_R(T-t)/\eta^2\), and letting \((t,y)\to(T,x)\) gives
\(v_*(T,x)\ge1=g_*(x)\).

(ii) The Hamiltonian \(H\) of \eqref{eq:Isaacs_Hamiltonian} inherits
from Assumption~\ref{ass:dynamics} the structure estimates
\[\begin{aligned}
    |H(t,x,p)-H(t,y,p)|&\le L_f\,|x-y|\,|p|,
    \\
    |H(t,x,p)-H(t,x,q)|&\le C_f\,(1+|x|)\,|p-q|,
\end{aligned}\]
where \(L_f\) is the Lipschitz constant and \(C_f\) the linear growth
constant of \(f\) (suprema and infima over the compact sets \(A\) and
\(B\) preserve a common modulus), and \(\sigma\) is Lipschitz in \(x\)
uniformly in \(t\).  These are the classical conditions under which
comparison holds for bounded semicontinuous sub- and supersolutions of
uniformly parabolic Bellman--Isaacs equations on the whole space; see
\cite[Thm.~8.2 and Sec.~5.D]{CrandallIshiiLions1992} for the parabolic
doubling argument and \cite{DaLioLey2006} for its implementation with
Lipschitz coefficients of linear growth on \([0,T]\times\R^d\).
Apply this argument in forward time \(s=T-t\), so the ordered terminal
traces become ordered initial traces. Quadratic spatial localization,
using \(\varepsilon e^{Cs}(1+|x|^2)\) with sufficiently large \(C\),
controls infinity under the linear-growth bounds. The initial ordering
excludes a positive maximum at \(s=0\); the interior viscosity inequalities
exclude one at \(s>0\). Letting \(\varepsilon\downarrow0\) gives comparison.

(iii) For bounded Lipschitz \(\psi\), the dynamic programming
principle \eqref{eq:DPP_general_payoff}, together with the standard
test-function argument and It\^o's formula, shows that \(V^\psi\) is a
viscosity solution of \eqref{eq:HJI_general_payoff} on
\([0,T)\times\R^d\); this is the setting of
\cite{FlemingSouganidis1989,BuckdahnLi2008}, and joint continuity of
\(V^\psi\) on \([0,T)\times\R^d\) follows as in \cite{BuckdahnLi2008}.
For bounded uniformly continuous \(\psi\), uniform Lipschitz approximation,
the bound \(\|V^\psi-V^\varphi\|_\infty\le\|\psi-\varphi\|_\infty\)
from Lemma~\ref{lem:value_monotonicity}, and viscosity stability give the
same conclusions.
Terminal attainment is again a consequence of
Lemma~\ref{lem:short_time_moment}: if \(\omega_\psi\) denotes a modulus
of continuity of \(\psi\) and \(|x|\le R\), then for every admissible
pair \((a,\gamma)\) and every \(\rho>0\),
\[\begin{aligned}
    \bigl|\E[\psi(X_T^{t,x;a,\gamma[a]})]-\psi(x)\bigr|
    &\le
    \E\bigl[\omega_\psi(|X_T^{t,x;a,\gamma[a]}-x|)\bigr]
    \\
    &\le
    \omega_\psi(\rho)
    +
    \frac{2\|\psi\|_\infty\,C_R\,(T-t)}{\rho^2}.
\end{aligned}\]
The bound is uniform in \((a,\gamma)\), hence survives the infimum over
\(a\) and the supremum over \(\gamma\); letting first \(t\uparrow T\) and
then \(\rho\downarrow0\) gives the locally uniform terminal attainment.
Uniqueness in the bounded class follows from part (ii) applied in both
directions.
\end{proof}

\subsection{Proof of the envelope error bound}

The error is controlled by the probability of reaching the target's
boundary tube at the terminal time.

\begin{proof}
Since \(g^+_\eps-g^-_\eps\) vanishes outside
\(\mathcal N_\eps(\partial\mathcal G)\) and takes values in \([0,1]\),
Assumption~\ref{ass:tube_estimate} gives, uniformly over admissible
pairs \((a,\gamma)\) and \((t,x)\in[0,T-\tau]\times K\),
\[
    0
    \le
    \E\!\left[
        g^+_\eps(X_T^{t,x;a,\gamma[a]})
        -
        g^-_\eps(X_T^{t,x;a,\gamma[a]})
    \right]
    \le
    \Prob\!\left(
        X_T^{t,x;a,\gamma[a]}
        \in
        \mathcal N_\eps(\partial\mathcal G)
    \right)
    \le\omega_{K,\tau}(\eps).
\]
Taking the infimum over \(a\) and then the supremum over \(\gamma\)
preserves this uniform additive bound, giving
\(0\le v^+_\eps-v^-_\eps\le\omega_{K,\tau}(\eps)\).
\end{proof}

\section{Policy iteration proofs}\label{app:howard}

The policy bounds follow from Hamiltonian ordering; convergence and
stagewise error estimates then use the residuals of the exact iterates.

\subsection{Policy bounds and fixed-pair gaps}\label{subsec:one_sided_policy_hamiltonians}
Use \(H_\beta\) and \(H^\alpha\) from Section~\ref{sec:PINNPI}.
Let \(v^-_{\eps,\beta}\) be the viscosity solution of the lower
one-sided policy equation
\begin{equation}
\label{eq:lower_policy_value}
\begin{cases}
    \partial_t u
    +
    H_\beta(t,x,\nabla_x u)
    +
    \dfrac12\Tr(\Sigma D^2_{xx}u)
    =0,
    & (t,x)\in[0,T)\times\R^d,\\[1mm]
    u(T,x)=g^-_\eps(x),
    & x\in\R^d.
\end{cases}
\end{equation}
Similarly, let \(v^{+,\alpha}_{\eps}\) be the viscosity solution of the
upper one-sided policy equation
\begin{equation}
\label{eq:upper_policy_value}
\begin{cases}
    \partial_t u
    +
    H^\alpha(t,x,\nabla_x u)
    +
    \dfrac12\Tr(\Sigma D^2_{xx}u)
    =0,
    & (t,x)\in[0,T)\times\R^d,\\[1mm]
    u(T,x)=g^+_\eps(x),
    & x\in\R^d.
\end{cases}
\end{equation}

\begin{proposition}
\label{prop:policy_wise_bracketing}
For disturbance and control feedbacks \(\beta\) and \(\alpha\) whose
policy equations admit bounded continuous viscosity solutions attaining
the stated terminal data,
\begin{equation}
\label{eq:policy_bracket_value}
    v^-_{\eps,\beta}(t,x)
    \le
    v^-_\eps(t,x)
    \le
    v(t,x)
    \le
    v^+_\eps(t,x)
    \le
    v^{+,\alpha}_{\eps}(t,x)
\end{equation}
for all \((t,x)\in[0,T]\times\R^d\).
\end{proposition}

\begin{proof}
If a smooth test function touches \(v^-_{\eps,\beta}\) from above,
its policy-operator value is nonnegative. Since \(H_\beta\le H\),
its full HJI-operator value is also nonnegative. Thus
\(v^-_{\eps,\beta}\) is a viscosity subsolution of the envelope equation.
The analogous test from below, using \(H\le H^\alpha\), makes
\(v^{+,\alpha}_\eps\) a viscosity supersolution. Comparison with the
matching continuous terminal data gives the outer inequalities;
\eqref{eq:value_bracketing} gives the middle ones.
\end{proof}

If both the control and disturbance policies are frozen, the
policy-evaluation equation becomes linear.  For a fixed pair
\((\alpha,\beta)\), define
\[
    q^{\alpha,\beta}(t,x,p)
    :=
    f(t,x,\alpha(t,x),\beta(t,x))\cdot p
\]
and
\[
    \mathcal H^{\alpha,\beta}[u]
    :=
    \partial_t u
    +
    q^{\alpha,\beta}(t,x,\nabla_x u)
    +
    \frac12\Tr(\Sigma D^2_{xx}u).
\]
A fixed-pair policy evaluation solves
\[
    \mathcal H^{\alpha,\beta}[u]=0.
\]
The fixed-pair Hamiltonian \(q^{\alpha,\beta}\) has no prescribed
ordering relative to \(H\).
Indeed, for an exact fixed-pair solution \(u\),
\[
    \mathcal H[u]
    =
    H(t,x,\nabla_x u)
    -
    q^{\alpha,\beta}(t,x,\nabla_x u).
\]
Thus a fixed-pair linear evaluation is not automatically a lower or
upper certificate.

On a verification cylinder \(Q\subset[0,T)\times\R^d\), define
\[
    \Delta^-_{\alpha,\beta,Q}[u]
    :=
    \sup_{(t,x)\in Q}
    \left(
        q^{\alpha,\beta}(t,x,\nabla_x u)
        -
        H(t,x,\nabla_x u)
    \right)_+,
\]
and
\[
    \Delta^+_{\alpha,\beta,Q}[u]
    :=
    \sup_{(t,x)\in Q}
    \left(
        H(t,x,\nabla_x u)
        -
        q^{\alpha,\beta}(t,x,\nabla_x u)
    \right)_+ .
\]
Assume these suprema are finite for \(Q=[0,T)\times\R^d\). They quantify
the correction from fixed-pair evaluation to the full HJI residual.

\subsection{Proof of exact Howard convergence}

Compactness allows the greedy policy identities to pass to the limit
in the HJI equation.

\begin{proof}
Proposition~\ref{prop:monotone_policy_certificates} gives
\(L_n\uparrow L_\infty\le v^-_\eps\) and
\(U_n\downarrow U_\infty\ge v^+_\eps\). By
Assumption~\ref{ass:exact_howard_regularity}, each subsequence has a further
subsequence converging in \(C^{1,2}\) on compact subcylinders. Uniqueness of
the pointwise limits therefore gives convergence of the whole sequences
in this topology.

The lower evaluation equation and greedy identity give
\[
\begin{aligned}
\mathcal H[L_{n+1}]
&=H(t,x,\nabla L_{n+1})-H_{\beta_{n+1}}(t,x,\nabla L_{n+1})\\
&=H(t,x,\nabla L_{n+1})-H(t,x,\nabla L_n)\\
&\quad+H_{\beta_{n+1}}(t,x,\nabla L_n)
-H_{\beta_{n+1}}(t,x,\nabla L_{n+1}).
\end{aligned}
\]
Both Hamiltonians are locally Lipschitz in the gradient, uniformly in the
policy, so the absolute value is bounded by
\(2C_K|\nabla L_{n+1}-\nabla L_n|\to0\) on each compact subcylinder.
Passing to the limit gives \(\mathcal H[L_\infty]=0\). Replacing
\(L_n,H_{\beta_{n+1}}\) by \(U_n,H^{\alpha_{n+1}}\) gives
\(\mathcal H[U_\infty]=0\).

The terminal values follow from
\[
 L_0\le L_\infty\le v^-_\eps,
 \qquad v^+_\eps\le U_\infty\le U_0.
\]
The bounds on each side are continuous up to \(T\) with the same terminal
data, by Assumption~\ref{ass:exact_howard_regularity} and
Theorem~\ref{thm:relaxed_wellposedness}(iii). Squeezing gives
\(L_\infty(T,\cdot)=g^-_\eps\) and
\(U_\infty(T,\cdot)=g^+_\eps\) continuously. Uniqueness in
Theorem~\ref{thm:relaxed_wellposedness}(ii)--(iii) identifies the limits as
\(v^-_\eps\) and \(v^+_\eps\). The bounds on \(v\) and the limit
\(\eps\downarrow0\) follow from \eqref{eq:value_bracketing}
and Proposition~\ref{prop:envelope_error}.
\end{proof}

\subsection{Proof of the stagewise pointwise certificate}
Residual corrections bracket each envelope value, after which the
envelope gap gives the error relative to the reachability value.

\begin{proof}
The inequalities \(0\le\mathcal H[L_n]\le r_n^-\) show that \(L_n\)
and \(L_n+r_n^-(T-t)\) are respectively a subsolution and a supersolution
with terminal value \(g^-_\eps\). Comparison brackets \(v^-_\eps\)
between them. Likewise, \(U_n-r_n^+(T-t)\le v^+_\eps\le U_n\).
The bounds relative to \(v\) follow from terminal-envelope bracketing
and, for the estimates involving \(\omega_{K,\tau}\),
Assumption~\ref{ass:tube_estimate} and Proposition~\ref{prop:envelope_error}.
\end{proof}

\section{Reference details}\label{app:benchmarks}

The feedback laws, slice processing, and reference coefficients below
specify the numerical comparisons in Section~\ref{sec:experiments}.

\subsection{Air3D feedbacks}

The rollouts use the explicit feedback laws
\[
    \omega_p^*(t,x)
    =
    \omega_{\max}\,\mathrm{sign}(p_{x_3}),
\]
and
\[
    \omega_e^*(t,x)
    =
    -\omega_{\max}\,
    \mathrm{sign}
    \left(
        x_2p_{x_1}-x_1p_{x_2}-p_{x_3}
    \right).
\]

\begin{figure}[htbp]
    \centering
    \begin{subfigure}[t]{0.32\columnwidth}
        \centering
        \includegraphics[width=\linewidth]{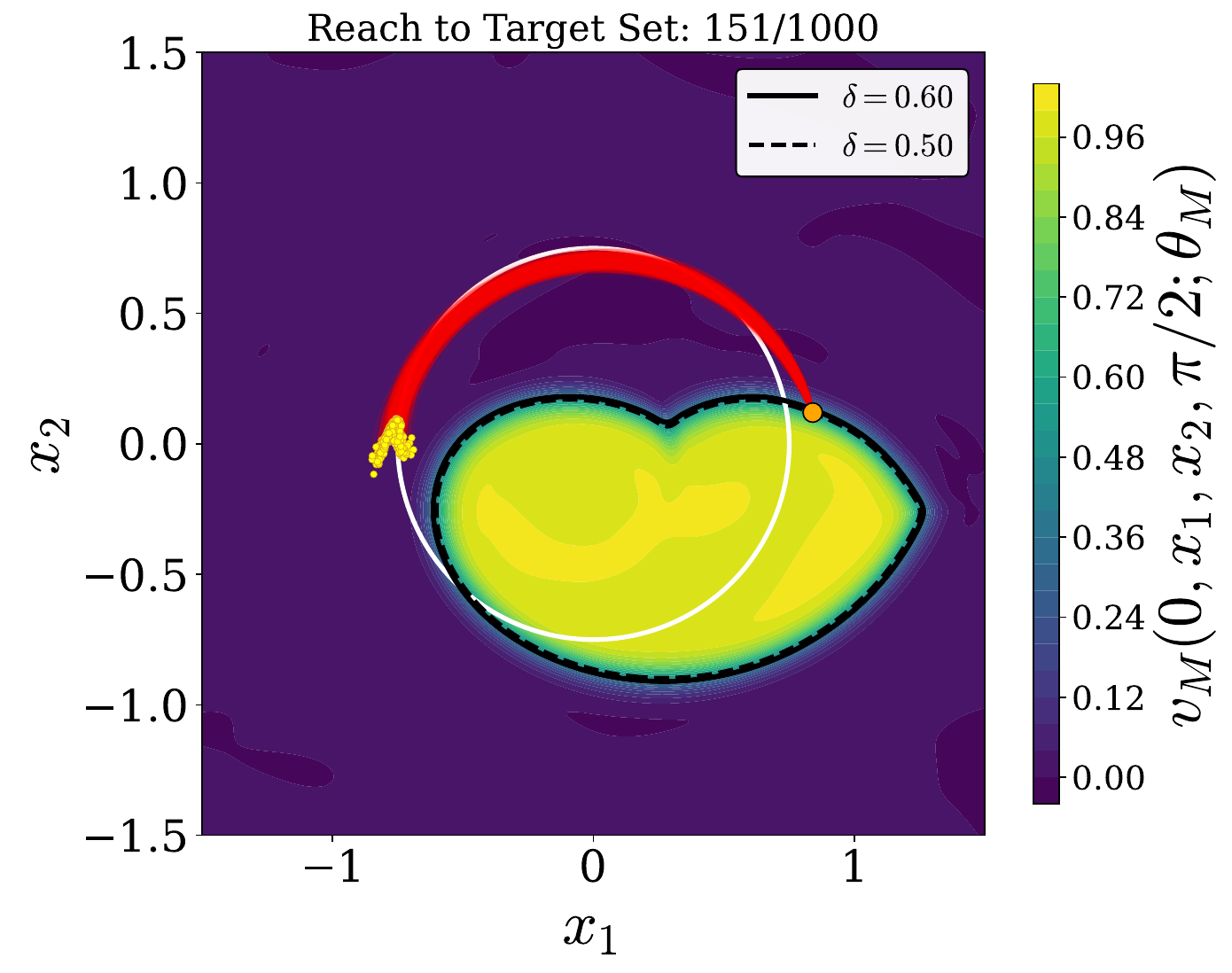}
        \caption{\(\kappa=0.5\)}
        \label{fig:air3d_mc_kappa05}
    \end{subfigure}\hfill
    \begin{subfigure}[t]{0.32\columnwidth}
        \centering
        \includegraphics[width=\linewidth]{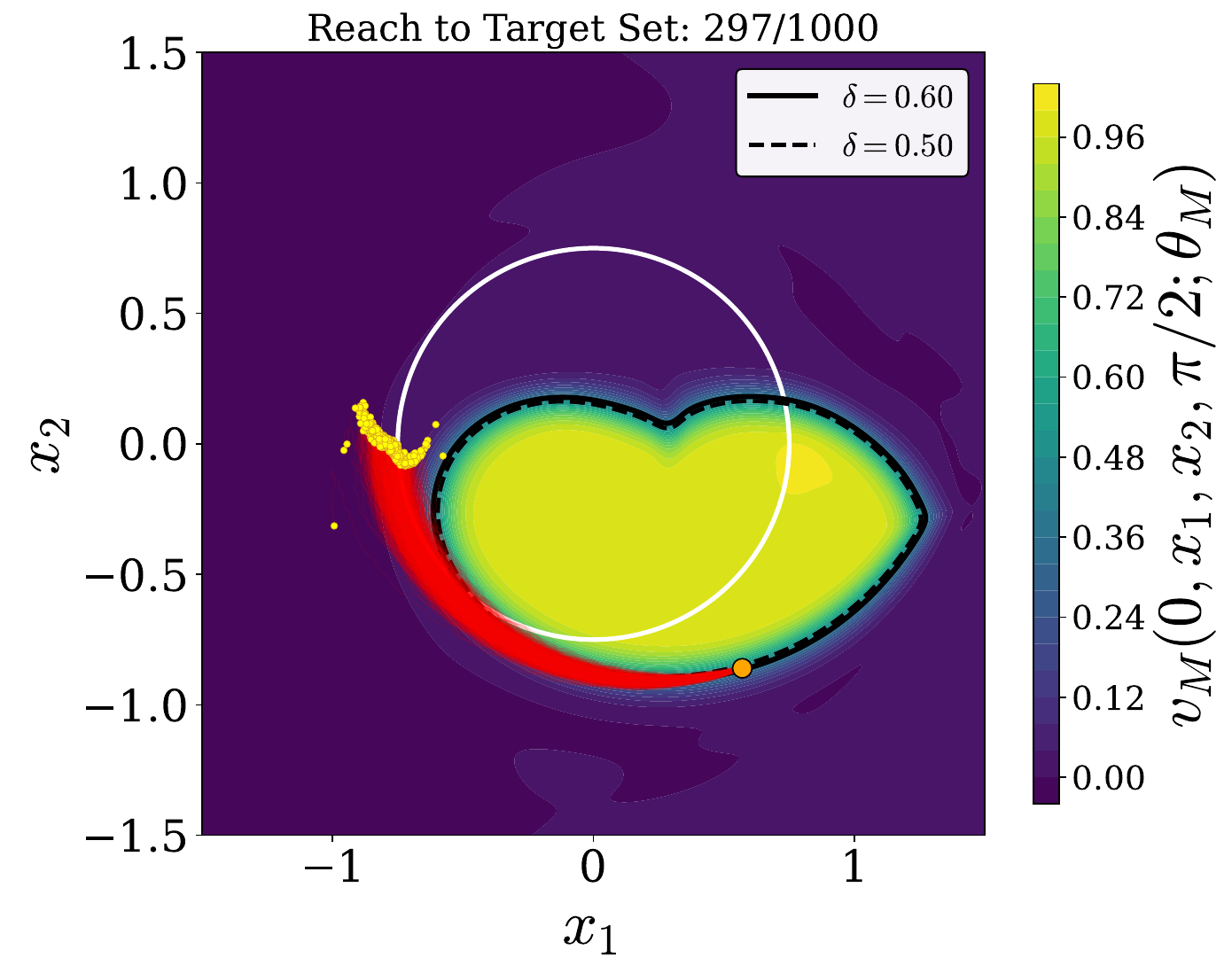}
        \caption{\(\kappa=1.0\)}
        \label{fig:air3d_mc_kappa10}
    \end{subfigure}\hfill
    \begin{subfigure}[t]{0.32\columnwidth}
        \centering
        \includegraphics[width=\linewidth]{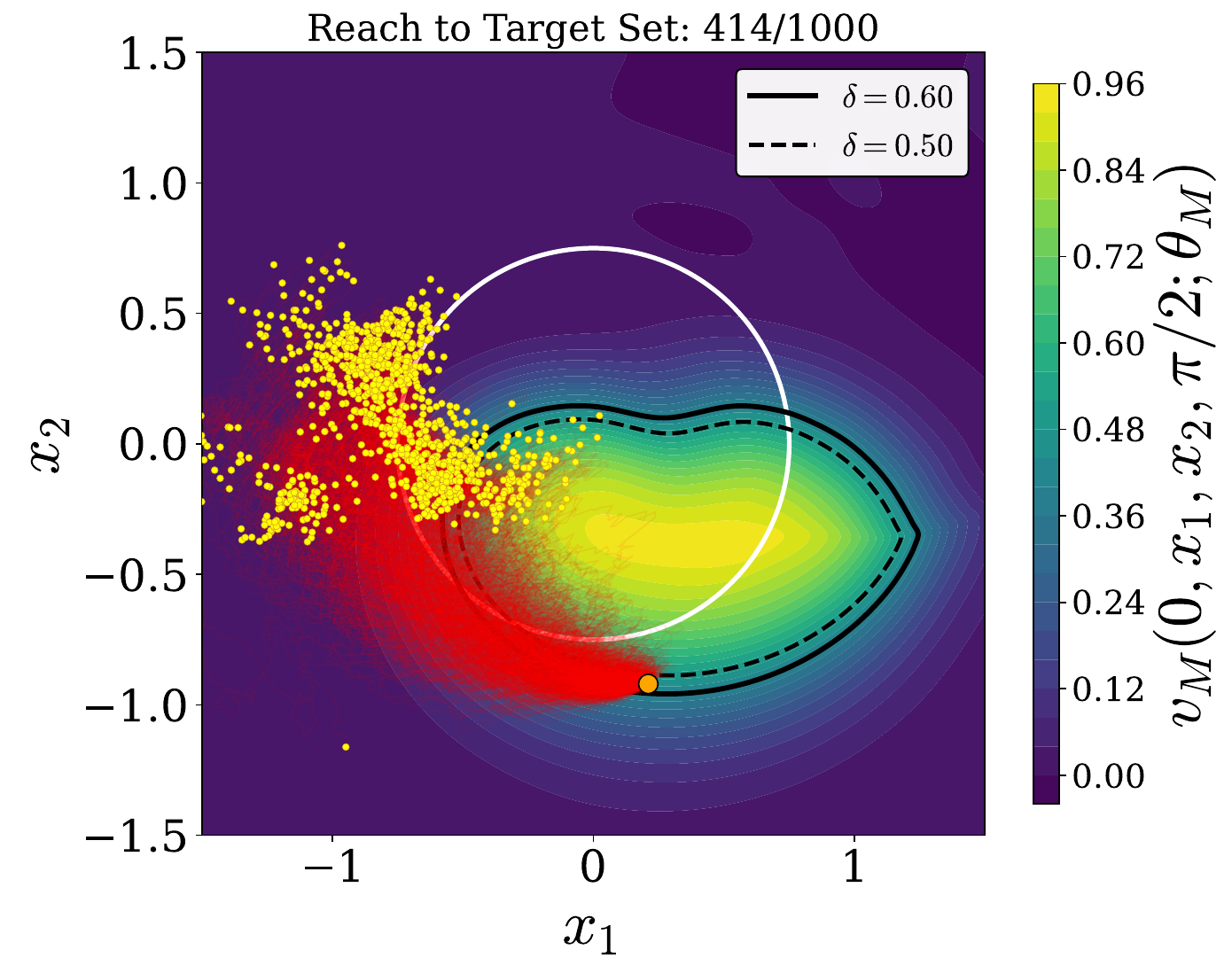}
        \caption{\(\kappa=5.0\)}
        \label{fig:air3d_mc_kappa50}
    \end{subfigure}
    \caption{Representative Monte Carlo rollouts for the learned value
    band \(v_\kappa(0,x)\in[0.4,0.5]\).  For each \(\kappa\), the plotted
    state is chosen so that its estimated success rate is closest to the
    empirical mean over the sampled states.  The solid and dashed
    contours indicate \(v_\kappa(0,x)=0.4\) and \(v_\kappa(0,x)=0.5\),
    respectively.}
    \label{fig:air3d_mc_kappa}
\end{figure}

Here \(p_{x_i}=\partial_{x_i}v_\kappa(t,x)\) is evaluated from the
learned candidate at the current time and state.
Figure~\ref{fig:air3d_mc_kappa} shows initial states whose empirical terminal
success rates are nearest to their scale-specific means.

\subsection{4D reference processing}

Before subsampling each 4D position slice, retain the connected component
containing the maximum at threshold 0.01, dilate its mask by eight fine-grid
pixels, and set values outside the mask to zero.

\begin{figure}[htbp]
    \centering
    \includegraphics[width=0.98\columnwidth]{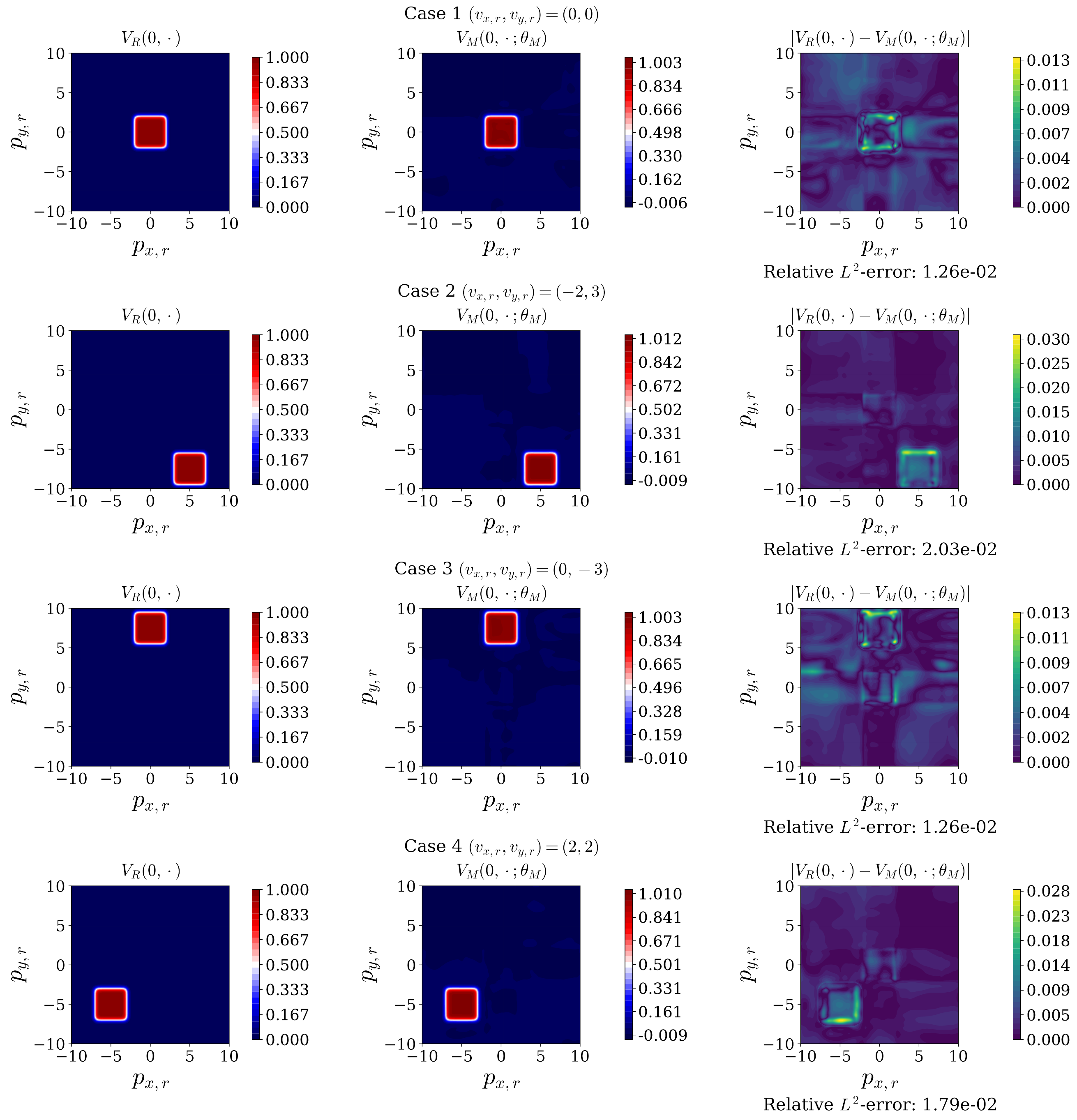}
    \caption{Reduced 4D quadrotor: value-function slices at \(t=0\).
    Each row shows, from left to right, the processed FDM reference \(v_R(0,\cdot)\),
    the PINN--PI result \(v_M(0,\cdot;\theta_M)\), and the absolute error
    \(|v_R(0,\cdot)-v_M(0,\cdot;\theta_M)|\) on the
    \((p_{x,r},p_{y,r})\)-plane for a fixed
    \((v_{x,r},v_{y,r})\).}
    \label{fig:quad4d_t0_cases_stacked}
\end{figure}

Figure~\ref{fig:quad4d_t0_cases_stacked} shows the value slices and
absolute errors for the four velocity pairs in Section~\ref{subsec:quad4d}.

\subsection{10D reference model}

The 10D comparison uses three subsystem references whose computational
domains and diffusion coefficients are specified below.

\begin{figure}[htbp]
    \centering
    \begin{subfigure}[t]{0.66\columnwidth}
        \centering
        \includegraphics[width=\linewidth]{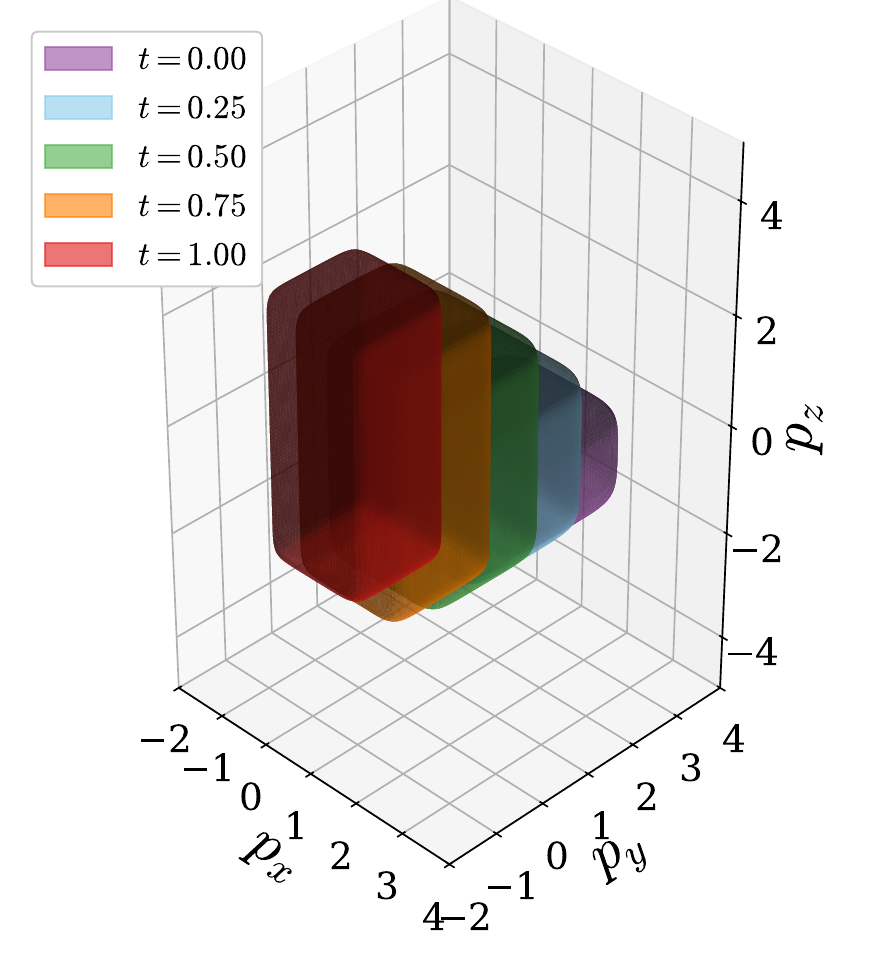}
    \end{subfigure}\hfill
    \begin{subfigure}[t]{0.33\columnwidth}
        \centering
        \includegraphics[width=\linewidth]{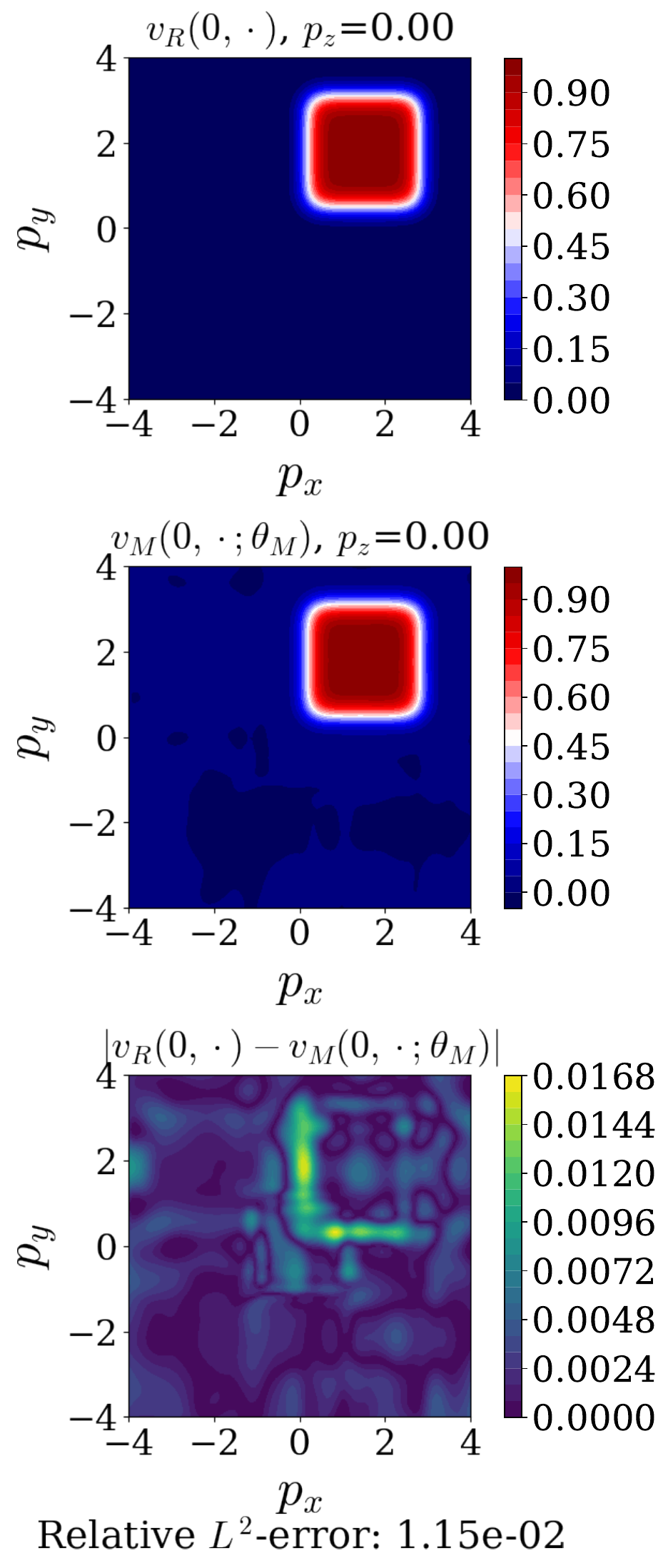}
    \end{subfigure}
    \caption{10D near-hover quadrotor example.
    Left: 3D isosurfaces of the PINN--PI value function
    \(v_M(t,x;\theta_M)\) in \((p_x,p_y,p_z)\) across time, with the
    remaining states fixed.  Right: a 2D slice in the \((p_x,p_y)\)-plane
    at \(p_z=0\), showing the FDM reference, the PINN--PI result, and the
    absolute error.}
    \label{fig:nearhover10d_slices}
\end{figure}

The computational domain is
\[
\begin{aligned}
(p_x,v_x,\theta_x,\omega_x)
&\in[-4,4]\times[-2,2]\times[-0.2,0.2]\times[-1,1],\\
(p_y,v_y,\theta_y,\omega_y)
&\in[-4,4]\times[-2,2]\times[-0.2,0.2]\times[-1,1],\\
(p_z,v_z)
&\in[-4,4]\times[-3,3].
\end{aligned}
\]
We use
\[
\sigma_x=\sigma_y=
\begin{bmatrix}
    0.12 & 0.03 & 0.01 & 0.02\\
    0.03 & 0.15 & 0.02 & 0.01\\
    0.01 & 0.02 & 0.10 & 0.01\\
    0.02 & 0.01 & 0.01 & 0.08
\end{bmatrix},
\qquad
\sigma_z=
\begin{bmatrix}
    0.15 & 0.02\\
    0.02 & 0.10
\end{bmatrix}.
\]

The restricted-product reference of Section~\ref{subsec:quad10d} uses the following blocks.

\[
x^{(1)}=(p_x,v_x,\theta_x,\omega_x),
\qquad
x^{(2)}=(p_y,v_y,\theta_y,\omega_y),
\]
\[
x^{(3)}=(p_z,v_z).
\]
Let \(p_{x^{(k)}}=\nabla_{x^{(k)}}v\) denote the gradient variable for
each subsystem.  With
\[
    |d_x|\le\bar d_x,\quad
    |d_y|\le\bar d_y,\quad
    |d_z|\le\bar d_z,
\]
\[
    |S_x|\le S_{x,\max},
    \qquad
    |S_y|\le S_{y,\max},
    \qquad
    T_{z,\min}\le T_z\le T_{z,\max},
\]
the subsystem Hamiltonians
\(H_k=\sup_{b^{(k)}}\inf_{a^{(k)}}f_k\cdot p^{(k)}\) have the closed
forms
\[
\begin{aligned}
H_1(x^{(1)},p)
&=
v_xp_{p_x}
+\bar d_x|p_{p_x}|
+g_{\rm grav} \tan(\theta_x)p_{v_x}
+(-d_1\theta_x+\omega_x)p_{\theta_x}
\\[-1mm]
&\hspace{3em}
-d_0\theta_xp_{\omega_x}
-n_0S_{x,\max}|p_{\omega_x}|,\\[1mm]
H_2(x^{(2)},p)
&=
v_yp_{p_y}
+\bar d_y|p_{p_y}|
+g_{\rm grav} \tan(\theta_y)p_{v_y}
+(-d_1\theta_y+\omega_y)p_{\theta_y}
\\[-1mm]
&\hspace{3em}
-d_0\theta_yp_{\omega_y}
-n_0S_{y,\max}|p_{\omega_y}|,\\[1mm]
H_3(x^{(3)},p)
&=
v_zp_{p_z}
+\bar d_z|p_{p_z}|
-g_{\rm grav} p_{v_z}
+k_TT_z^\ast p_{v_z},
\end{aligned}
\]
where
\[
    T_z^\ast=T_{z,\max}\quad\text{if }p_{v_z}\le0,
    \qquad
    T_z^\ast=T_{z,\min}\quad\text{if }p_{v_z}>0,
\]
so that \(k_TT_z^\ast p_{v_z}
=\inf_{T_{z,\min}\le T_z\le T_{z,\max}}k_TT_zp_{v_z}\).  In accordance
with the convention \eqref{eq:Isaacs_Hamiltonian}, the disturbance
bounds enter with positive sign and the control bounds with negative
sign.  (The level-set reference computation uses the time-to-go
variable \(s=T-t\) and therefore the negated Hamiltonians \(-H_k\).)
The full Hamiltonian is \(H=H_1+H_2+H_3\), and the diffusion operator is
treated blockwise for the reference computation.

Figure~\ref{fig:nearhover10d_slices} shows additional isosurfaces and the reference slice specified in Section~\ref{subsec:quad10d}.

\section{Factorization for product stochastic games}
\label{app:decomp}

The reference construction uses a restricted game with product dynamics,
independent Brownian blocks, a product target, and blockwise admissible
strategies.

\subsection{Product structure}

Partition \(x=(x^{(1)},\dots,x^{(K)})\), with
\(x^{(k)}\in\R^{d_k}\) and \(\sum_kd_k=d\). The action sets are
\(A=\prod_kA_k\) and \(B=\prod_kB_k\), with components
\(a^{(k)}\) and \(b^{(k)}\).

\begin{assumption}
\label{ass:app_decomp_dyn}
For each \(k=1,\dots,K\), the \(k\)-th subsystem evolves according to
\[
    dX^{(k)}_s
    =
    f_k\!\left(
        s,
        X^{(k)}_s,
        a^{(k)}_s,
        b^{(k)}_s
    \right)\,ds
    +
    \sigma_k\!\left(
        s,
        X^{(k)}_s
    \right)\,dW^{(k)}_s,
\]
where \(W^{(k)}\) is an \(m_k\)-dimensional Brownian motion.  The
Brownian motions \(W^{(1)},\dots,W^{(K)}\) are mutually independent, and
\[
    \sum_{k=1}^K m_k=m.
\]
Equivalently, the full diffusion matrix is block diagonal:
\[
    \sigma(s,x)
    =
    \diag\!\left(
        \sigma_1(s,x^{(1)}),
        \dots,
        \sigma_K(s,x^{(K)})
    \right).
\]
\end{assumption}

\begin{assumption}
\label{ass:app_decomp_target}
There exist closed sets \(\mathcal G_k\subset\R^{d_k}\) such that
\[
    \mathcal G=\prod_{k=1}^K \mathcal G_k .
\]
Equivalently,
\[
    g(x)
    =
    \1_{\mathcal G}(x)
    =
    \prod_{k=1}^K g_k(x^{(k)}),
    \qquad
    g_k:=\1_{\mathcal G_k}.
\]
\end{assumption}

Admissible controls and strategies are blockwise:
\(a=(a^{(1)},\dots,a^{(K)})\) and
\(\gamma[a]=(\gamma_1[a^{(1)}],\dots,\gamma_K[a^{(K)}])\).
Each \(\gamma_k\) is nonanticipative, and each subsystem control and
strategy uses only its own Brownian filtration, with no shared
randomization. Define the restricted game and subsystem values by
\[\begin{aligned}
 V^{\rm prod}(t,x)&:=\sup_{(\gamma_k)}\inf_{(a^{(k)})}
       \E\!\left[\prod_{k=1}^K g_k(X_T^{(k)})\right],\\
 v_k(t,x^{(k)})&:=\sup_{\gamma_k}\inf_{a^{(k)}}
       \E[g_k(X_T^{(k)})].
\end{aligned}\]

\subsection{Value factorization}

Independence separates the expected terminal payoff; the following
identity separates the optimization of its nonnegative factors.

\begin{lemma}
\label{lem:app_sep_extrema}
Let \(\mathcal U_k\) and \(\mathcal V_k\) be nonempty sets, and let
\[
    F_k:\mathcal U_k\times\mathcal V_k\to[0,\infty)
\]
be bounded nonnegative functions.  Then
\[
    \sup_{v_1,\dots,v_K}
    \inf_{u_1,\dots,u_K}
    \prod_{k=1}^K F_k(u_k,v_k)
    =
    \prod_{k=1}^K
    \left(
        \sup_{v_k}
        \inf_{u_k}
        F_k(u_k,v_k)
    \right).
\]
\end{lemma}

\begin{proof}
For a fixed tuple \((v_k)\), independent minimization of the nonnegative
factors gives \(\inf_{(u_k)}\prod_kF_k(u_k,v_k)
=\prod_k\inf_{u_k}F_k(u_k,v_k)\), using approximating minimizers if needed.
Independent maximization of these nonnegative factors gives the stated
identity, again by approximation when extrema are not attained.
\end{proof}

\begin{proposition}
\label{prop:app_factorization}
Under Assumptions~\ref{ass:app_decomp_dyn}--\ref{ass:app_decomp_target}
and under the restriction to product controls and product disturbance
strategies,
\[
    V^{\rm prod}(t,x)
    =
    \prod_{k=1}^K
    v_k(t,x^{(k)}),
    \qquad
    (t,x)\in[0,T]\times\R^d.
\]
\end{proposition}

\begin{proof}
For each fixed admissible product control and strategy, the subsystem
processes are independent: their dynamics are uncoupled and each block
uses its own independent Brownian motion without shared randomization.
Thus, with \(F_k(a^{(k)},\gamma_k):=\E[g_k(X_T^{(k)})]\),
\[
 \E[g(X_T)]=\prod_{k=1}^K F_k(a^{(k)},\gamma_k).
\]
Since \(0\le F_k\le1\), Lemma~\ref{lem:app_sep_extrema} separates the
infimum and supremum, giving \(V^{\rm prod}=\prod_kv_k\).
\end{proof}

An unrestricted disturbance strategy can couple blocks through the
observed control history, so this identity concerns the product game.
The proof also applies to bounded nonnegative terminal factors, including
the smoothed payoffs used in the experiments.

\subsection{Hamiltonian factorization}

For product action sets, the block drift gives
\(f(t,x,a,b)\cdot p=\sum_k f_k(t,x^{(k)},a^{(k)},b^{(k)})\cdot p^{(k)}\).
Optimizing each block separately yields
\[\begin{aligned}
 H(t,x,p)&=\sum_{k=1}^K H_k(t,x^{(k)},p^{(k)}),\\
 H_k(t,x^{(k)},p^{(k)})&:=\sup_{b^{(k)}\in B_k}\inf_{a^{(k)}\in A_k}
 f_k(t,x^{(k)},a^{(k)},b^{(k)})\cdot p^{(k)}.
\end{aligned}\]

\begin{lemma}
\label{lem:subsystem_isaacs}
Suppose the product drift and product action-set assumptions above hold.
If the full Hamiltonian satisfies the Isaacs condition
\[
    \sup_{b\in B}
    \inf_{a\in A}
    f(t,x,a,b)\cdot p
    =
    \inf_{a\in A}
    \sup_{b\in B}
    f(t,x,a,b)\cdot p
\]
for every \((t,x,p)\), then each subsystem Hamiltonian satisfies
\[
    \sup_{b^{(k)}\in B_k}
    \inf_{a^{(k)}\in A_k}
    f_k(t,x^{(k)},a^{(k)},b^{(k)})\cdot p^{(k)}
    =
    \inf_{a^{(k)}\in A_k}
    \sup_{b^{(k)}\in B_k}
    f_k(t,x^{(k)},a^{(k)},b^{(k)})\cdot p^{(k)}.
\]
\end{lemma}

\begin{proof}
Fix \(k\) and choose a full gradient vector \(p\) such that
\(p^{(j)}=0\) for \(j\ne k\).  Then
\[
    f(t,x,a,b)\cdot p
    =
    f_k(t,x^{(k)},a^{(k)},b^{(k)})\cdot p^{(k)}.
\]
Applying the full Isaacs identity to this choice of \(p\) gives exactly
the subsystem Isaacs identity.
\end{proof}

\section{Computational audit}\label{app:computational_audit}

\subsection{Scope and provenance}
The audit uses the public repository
\url{https://github.com/HeeJunYang09/certified-probabilistic-reachability},
commit \nolinkurl{0b9436a2f156ce9ce9509f26d7972fff2f59b983}.
All 158 release-manifest files pass their checksum checks. The ten figure PDFs used in this manuscript are byte-identical to the published aliases in this repository. The external 4D reference was downloaded and its documented SHA-256 verified. Scientific source files and archived checkpoints were not edited.

The execution environment is Windows, CPU, Python 3.13.15, JAX/jaxlib 0.11.0, NumPy 2.4.6, and python-flint 0.9.0. The archive specifies python-flint 0.6.0, which could not be installed as a Windows wheel in this environment. A compatibility loader removes only obsolete, empty \texttt{named\_shape} metadata from older JAX pickles; it refuses nonempty metadata. Canonical 1D parameter hashes are unchanged. Air3D path reproduction uses \texttt{jax\_threefry\_partitionable=False}, the historical PRNG convention.

\subsection{Figure and table correspondence}
\begin{longtable}{@{}p{.16\linewidth}p{.79\linewidth}@{}}
\caption{Correspondence between manuscript results and computational checks.}
\label{tab:audit_figure_checks}\\
\toprule Item & Executed verification \\ \midrule
\endfirsthead
\toprule Item & Executed verification \\ \midrule
\endhead
\bottomrule
\endlastfoot
Fig.~\ref{fig:1d_certified_bracket} & Frozen 1D candidates, stored residual and boundary evidence, and FDM data regenerate the certified-bracket figure. The two finest FDM solves were also recomputed. \\
Fig.~\ref{fig:air3d_pinn_fdm_error} & The Air3D checkpoint and archived finite-difference array regenerate the value and error panels. \\
Figs.~\ref{fig:quad4d:t0:level05:slices} and~\ref{fig:quad4d_t0_cases_stacked} & The 4D checkpoint and downloaded reference regenerate the contour and error panels, including the documented reference post-processing. \\
Fig.~\ref{fig:nearhover10d_slices} & The 10D adaptive checkpoint regenerates the isosurfaces and the comparison on the archived reference slice. \\
Fig.~\ref{fig:air3d_mc_kappa} & Stored trajectories redraw the panels; their terminal success counts were recalculated directly. Fresh simulation checks are reported below. \\
Fig.~\ref{fig:nh10d_adaptive_random_t0_t05_t10} & Adaptive and uniform 10D checkpoints regenerate the sampling-comparison panels. \\
Table~\ref{tab:nh10d_l2_time_adaptive_random} & The evaluation script recomputes all five reference times for both checkpoints. \\
\end{longtable}
The figure regeneration checks establish correspondence of code, frozen candidates, and archived data. They do not establish that every checkpoint can be reproduced by a fresh training run.

\subsection{Scalar certificate and independent reference}
The stored evidence audit examines all 202,081 lower and 239,406 upper accepted cells. Exact rational checks show that each cover tiles the positive-half cylinder without gaps or overlaps. Canonical parameter hashes, stored residual upper bounds below 0.1, terminal construction, lateral records, and set-classification records pass.

A fresh outward-rounded Taylor/Arb computation evaluates 64 selected cells per side at 192-bit precision. The selection includes the eight largest stored violations, unresolved-gradient cells where available, endpoints of the record, and a seeded random subset. All 128 recomputed cells pass the 0.1 cap. Both complete lateral-boundary checks were recomputed, with lower correction 0 and upper correction 0.001. A separate fresh set-classification run uses the audited archived residual certificate as its prerequisite. This is not a fresh recomputation of every interior residual enclosure. The archived full residual computation and its mathematical enclosure implementation remain part of the trusted evidence.

The two monotone finite-difference solves at spacing 0.00125 were rerun. The maximum discrepancy from the archived solution arrays is $2.22\times10^{-16}$. These solves use terminal values $G_\pm=g_\varepsilon^\pm\pm\beta$ and boundary values $-1$ and $1$ on $(-3.3,3.3)$. Their gap includes boundary and terminal modifications. It must not be identified solely with the whole-space terminal-envelope smoothing gap.

\subsection{Four-dimensional reference cross-check}
For the implemented bounds $u_{\max}=d_{\max}=1$, the optimized control terms cancel. Each position--velocity block therefore has linear stochastic drift. If $h=T-t$ and $Q=\sigma_x\sigma_x^\top$, the final position has mean $p+hv$ and variance
\[s^2=hQ_{11}+h^2Q_{12}+h^3Q_{22}/3.\]
At $h=2.5$, $s^2=0.04985416666666667$. The Gaussian expectation of $\frac12[1+\tanh((2-|Z|)/0.2)]$, multiplied across the two independent blocks, gives an independent whole-space comparator. Quadrature orders 128 and 256 differ by less than $1.8\times10^{-11}$ in maximum norm on the evaluated slices. This convergence check is numerical, not interval certification.

The original plotting code retains the component containing the maximum above threshold 0.01, dilates it by eight four-neighbor grid steps, and zeros values outside it before subsampling the $1601^2$ grid to $201^2$. The table reports neural relative $L^2$ errors against each comparator.
\begin{center}\begin{tabular}{rrrr}\toprule
Velocity & Raw FDM & Processed FDM & Gaussian \\ \midrule
$(0,0)$ & 0.012547 & 0.012580 & 0.012550 \\
$(-2,3)$ & 0.020166 & 0.020287 & 0.020163 \\
$(0,-3)$ & 0.012526 & 0.012551 & 0.012528 \\
$(2,2)$ & 0.017752 & 0.017867 & 0.017751 \\
\bottomrule\end{tabular}\end{center}
Raw FDM relative discrepancies against the Gaussian comparator are below $8\times10^{-5}$ on these four slices. The processing changes neural relative errors by at most approximately $1.22\times10^{-4}$ in absolute terms. This parameter choice provides a useful reference check but does not test a nontrivial competition between 4D control and disturbance bounds.

\subsection{Air3D Monte Carlo reproduction}
The experiment fixes $x_3=\pi/2$ and selects initial positions from a $301\times301$ grid with learned value in $[0.4,0.5]$. It uses the terminal disk of radius 0.75, Euler--Maruyama time step 0.01, and 1,000 rollouts for each seed 0--999 at each noise scale. It evaluates learned feedbacks and an indicator payoff; the trained PDE uses a smoothed terminal payoff. The learned level therefore has no independently established calibration guarantee.

The saved arrays give means 0.150734, 0.297417, and 0.414437, with population standard deviations 0.085126, 0.088437, and 0.045277. All six saved representative trajectory files give their recorded terminal rates when counted directly. A fresh simulation under the current default JAX PRNG mode produces different paths; the historical mode recovers the representative success counts. Floating-point changes near a bang--bang switching surface can change individual trajectories, so bitwise path identity is not claimed.

The following table compares completed full initial-state scans. ``Changed'' counts initial states whose integer number of terminal successes differs from the archive; the last column is the largest absolute count difference among their 1,000 paths.
\begin{center}\begin{tabular}{rrrrr}\toprule
$\kappa$ & Archived mean & Fresh mean & Changed & Max. count difference \\ \midrule
0.5 & 0.150734 & 0.150731 & 3 & 1 \\
1.0 & 0.297417 & 0.297417 & 0 & 0 \\
5.0 & 0.414437 & 0.414437 & 0 & 0 \\
\bottomrule\end{tabular}\end{center}
All three scans completed: a total of 3,000,000 fresh rollouts.

\subsection{Ten-dimensional reference comparison and adaptive sampling}
At $t=0$, recomputed relative $L^2$ errors are 0.011451982 for adaptive collocation and 0.996288836 for uniform collocation. At $t=0.25,0.5,0.75$, the adaptive errors are 0.010111762, 0.009128571, and 0.008694162. The terminal discrepancy is approximately $1.04\times10^{-7}$ in this environment, compared with $1.09\times10^{-7}$ in the archived report; both candidates satisfy their hard terminal ansatz by construction. This terminal row measures numerical reference/rounding consistency, rather than learned interior accuracy.

The network accepts ten state coordinates, but the quantitative errors above are measured on two-dimensional slices with the other coordinates fixed. The appendix's stochastic factorization is proved for restricted product strategies; that proposition alone does not establish that the product reference is the unrestricted game value. Neither the slices nor the 3D isosurfaces constitute a ten-dimensional certificate or a generic-dimensional scaling benchmark.

Adaptive and uniform training smoke runs pass. Full training was not repeated: the archived 10D runs each report approximately 15.4 hours on their original accelerator. The 1D, Air3D, and 4D archived training schedules were likewise not rerun from initialization. The public 4D training source initializes seed 42, while the frozen checkpoint filename contains ``seed0''; filenames alone cannot resolve the complete provenance of that original training run.

The original MATLAB reference generators for the Air3D, 4D, and 10D data are not present in the public repository. Their archived arrays can be checked and used, and the 4D arrays have the independent Gaussian check above, but the entire original MATLAB workflow is not reproducible from this repository alone.

\section*{Funding}
This work was supported by the National Institute for Mathematical Sciences
(NIMS) grant funded by the Korea government (MSIT) (No. B26810000), and by the
National Research Foundation of Korea (NRF) grants funded by the Korea
government (MSIT) (RS-2023-00219980 and RS-2026-25490291).
Yeoneung Kim was also supported by Yonsei University.

\section*{Data and code availability}
\url{https://github.com/HeeJunYang09/certified-probabilistic-reachability}

\section*{Declaration of competing interest}
The authors declare that they have no known competing financial interests
or personal relationships that could have appeared to influence the work
reported in this paper.

\section*{CRediT authorship contribution statement}
\textbf{Minjung Gim:} Methodology, Software, Formal analysis, Investigation,
Writing -- original draft.
\textbf{Hee Jun Yang:} Methodology, Software, Formal analysis, Investigation,
Writing -- original draft.
\textbf{Yeoneung Kim:} Supervision, Funding acquisition,
Writing -- review \& editing.

\section*{Declaration of generative AI and AI-assisted technologies in the manuscript preparation process}
OpenAI Codex assisted with manuscript revision and preparation of
computational verification scripts. The authors take full responsibility
for the mathematical arguments, results, and final text.

\bibliographystyle{unsrtnat}
\bibliography{references}

\end{document}